\documentclass[12pt,graphicx]{amsart}

\usepackage{amscd,amssymb,amsthm,latexsym, amsfonts,graphicx,amsmath}
\newtheorem{thm}{Theorem}[section]

\newtheorem{lem}[thm]{Lemma}

\newtheorem{defn}{{\sc Definition}}[section]

\newtheorem{rem}{{\sc Remark}}[section]

\newcommand{\bas}{\begin{eqnarray*}}

\newcommand{\eas}{\end{eqnarray*}}

\newcommand{\ba}{\begin{eqnarray}}

\newcommand{\ea}{\end{eqnarray}}



\begin{document}




























\title{A Gluing Theorem for Special Lagrangian Submanifolds}
\author{Sema Salur}

\vspace{.1in}

\address{Department of Mathematics, Cornell University, Ithaca, NY 14850}

\email{salur@math.cornell.edu}

\maketitle
\scriptsize

{\bf Abstract.}  The purpose of this paper is to prove a gluing theorem for a given special Lagrangian submanifold of a Calabi-Yau 3-fold. The proof will be an adaption of the gluing techniques in J-holomorphic curve theory, [7], [10], [12], [15]. It is a well known procedure in geometric analysis to construct new solutions to a given nonlinear partial differential equation by gluing known solutions [14]. First an approximate solution is constructed and then using analytic methods it is perturbed to a real solution. In this paper the gluing theorem will be used for smoothing a singularity of a special Lagrangian submanifold. In particular, we will show that given a special Lagrangian submanifold $L$ of a Calabi-Yau manifold $X$ with a particular codimension-two self intersection $K$ it can be approximated by a sequence of smooth special Lagrangian submanifolds and therefore $L$ is a limit point in the moduli space.  

\vspace{.1in}

\normalsize

\section{Introduction}

Special Lagrangian submanifolds has become an important subject of study [1],[3],[4],[5],[11] due to recent developments in Mirror Symmetry. In 1996, Strominger, Yau and Zaslow [18] proposed a geometric construction of the mirror manifold via special Lagrangian tori fibration. They claimed that the mirror of a Calabi-Yau can be obtained by some suitable compactification of the dual of this fibration. Therefore to find a compactification and understand the close relations with the mirror symmetry one first should understand the singularities of the moduli space of special Lagrangian submanifolds and how to characterize them.

\vspace{.1in}

Another motivation to study this moduli space is the hope to obtain new invariants for symplectic manifolds. In [13], McLean showed that the moduli space of all infinitesimal deformations of a smooth, compact, orientable special Lagrangian submanifold $L$ in a Calabi-Yau manifold $X$ within the class of special Lagrangian submanifolds is a smooth manifold of dimension equal to $b_1(L)$, the first Betti number of $L$. So if we can obtain some compactification of this moduli space then we will be able to count rational homology special Lagrangian spheres (hoping that there are finitely many of them) in a Calabi-Yau manifold [8]. McLean's result has been extended to symplectic manifolds by the author [16] and with a suitable compactification we can also obtain invariants for symplectic manifolds. Unfortunately, this program is far from complete and so far we don't know much about the global picture. With this work we hope to give some idea about possible singularities that can occur in this moduli space. 

\vspace{.1in}

In this paper we will consider a $3$-dimensional compact special Lagrangian submanifold $L$ in a Calabi-Yau manifold $X^{6}$ and we will assume that $L$ has a codimension two irreducible self-intersection $K$. The local model for this intersection will be explained in detail in section 2. For simplicity we will work in an open ball $V$ around $K$ and call the two parts of $L$ which come together and intersect along $K$, as $L_1$ and $L_2$. We will also assume that $L_1$ and $L_2$ intersect perpendicularly along $K$ and the normal bundles have opposite Euler classes: $e(N_{L_1}K)+e(N_{L_2}K)=0$. Also note that by McLean's result [13] the corresponding linearized operators $D_{L_1}$ and $D_{L_2}$ are surjective and this means that the corresponding moduli spaces of infinitesimal deformations of $L_1$ and $L_2$ are smooth manifolds near $L_1$ and $L_2$, respectively. Under the given assumptions our goal will be then to show that $L$ is a limit point in the moduli space of special Lagrangian submanifolds. 

\vspace{.1in}

More precisely for a given singular special Lagrangian submanifold $L$, we will first construct an approximate special Lagrangian submanifold in an open ball $V$ around the singular set $K$ and use the Implicit Function Theorem to prove that there exists a {\em true} special Lagrangian submanifold nearby. Inside $V$ the smoothing of the singularity can be represented as ${L_1\#L_2}$. In order to prove the existence of a {\em true} special Lagrangian we need to get a uniform estimate for the right inverse for the linearized operator $D_{L_1\#L_2}$.

\vspace{.1in}

In this paper we will prove the following theorem for the simplest type of singularity of the form $z_1.\overline {z_2}=0$.

\vspace{.1in}

{\thm: Given a connected immersed special Lagrangian submanifold $L^3$ of a Calabi-Yau manifold $X^6$ with a particular irreducible self intersection $K$ of codimension-two (singularity of type $z_1.\overline {z_2}=0$) it can be approximated by a sequence of smooth special Lagrangian submanifolds and therefore $L$ is a limit point in the moduli space.}

\vspace{.1in}

{\rem :} We need the irreducibility condition only in proving the eigenvalue estimate in section 4. Other parts of the proof do not require this condition. We hope to extend our result to the reducible case later. Also our gluing construction can be generalized to higher dimensional special Lagrangians $L^n$, $(n>3)$ by taking $K$ to be of dimension $n-2$ but since there are some regularity problems to be resolved in higher dimensions [9], in this paper we did our construction only for Calabi-Yau 3-folds.

\vspace{.1in}

First we recall some basic definitions:

\vspace{.1in}

 {\defn : A Calabi-Yau manifold $X^n$ is a K\"{a}hler manifold of complex dimension $n$ with a covariant constant holomorphic $n$-form. (equivalently it is a Riemannian manifold with holonomy contained in SU(n))}.

\vspace{.2in}

One other equivalent definition for Calabi-Yau manifolds is that they are K\"{a}hler manifolds with first Chern class $c_1=0$.

\vspace{.2in}

Calabi-Yau manifolds are equipped with a K\"{a}hler 2-form $\omega$, an almost complex structure $J$ which is tamed by $\omega$, the compatible Riemannian metric $\tilde{g}$ and a non-vanishing holomorphic $(n,0)$-form $\xi$.

\vspace{.2in}

{\defn : An $n$-dimensional submanifold $L\subseteq X$ is special Lagrangian if $L$ is Lagrangian (i.e. $\omega|_L\equiv 0$) and $Im(\xi)$ restricts to zero on $L$. Equivalently, $Re(\xi)$ restricts to be the volume form on $L$ with respect to the induced metric.}[6]

\section{Connected Sums of Special Lagrangian Submanifolds}

In this section, we will first define the local model for the singularity and for a given gluing parameter $\delta$, we will define a gluing map, smooth the singularity, and construct a smooth approximate special Lagrangian submanifold $H_\delta = L_1\#_\delta  L_2$ nearby. Recall that inside an open ball $V$ around $K$, the two parts of $L$ which come together and intersect along $K$ can be studied as a pair of special Lagrangian submanifolds $L_1$ and $L_2$ with opposite Euler classes. Here our goal is to construct a one-parameter family of approximate special Lagrangian which agrees with $L_1$ and $L_2$ outside a tubular neighbourhood of their intersection $K$. Our model locally consists of two special Lagrangian immersions $l_1$ and $l_2$,

\vspace{.1in}

$l_i: L=S^1\times S^1\times K\rightarrow L_i \subset X$ $  i=1,2$

\vspace{.1in}

\noindent intersecting orthogonally along $K$. The local geometry of $K$ will be explained in detail later. The connected sum $L_1$$\#_{\delta}$$L_2$ is obtained by removing the tubular neighbourhoods of $K$ in $L_1$ and $L_2$, and joining the boundaries of these cylinders. It is clear that for $L_1\#_{\delta}L_2$ to have an orientation compatible with the given orientations of $L_1$ and $L_2$, the boundaries of the tubular neighbourhoods must be joined by an orientation reversing map. The condition that the special Lagrangian submanifolds have opposite Euler classes and their intersection $K$ is codimension-$2$ guarantees the existence of an orientation-reversing diffeomorphism.



\subsection {PreGluing}

\vspace{.1in}

We are given a pair of three dimensional special Lagrangian submanifolds $L_1=T_1^2\times K$ and  $L_2=T_2^2\times K$ which are intersecting along $K$. This intersection can be smoothed at one of its points ${\mathcal Y}= (\kappa,0,0)$, $\kappa\in K$ as in the figure (1). Here $T_1=l_1(S^1\times S^1)$ and $T_2=l_2(S^1\times S^1)$ both represent two dimensional torus.

\begin{figure}[h]

\includegraphics{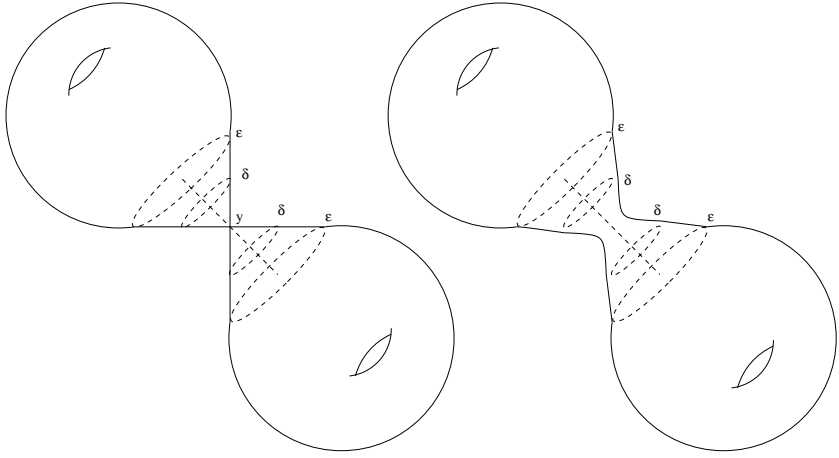}

\caption{\small Smoothing of $T_1\wedge T_2$ at one point}

\label{sevil-1}

\end{figure}

\vspace{.1in}

The construction of the approximate special Lagrangian submanifold is as follows:

\vspace{.1in}

Let $\mathcal D_1\subset T_1$ and $\mathcal D_2\subset T_2$ be two disks around ${\mathcal Y}$. Then a small neighbourhood around ${\mathcal Y}$ will be $(K\times \mathcal D_1\times \{0\}) \cup  (K\times \{0\}\times \mathcal D_2)$, where the singularity has the form  $K\times \{0\}\times \{0\}$. We can identify the small disks $\mathcal D_1$ and $\mathcal D_2$ with the complex plane {\bf C}.

\vspace{.1in}

Let $z_1=(x,y)$ and $z_2=(u,v)$ be the complex valued coordinates on $T_1$ and $T_2$, respectively. Then $T_1\wedge T_2$ is given by identifying $z_1=0$ in $T_1$ with $z_2=0$ in $T_2$. We then cut off $|z_1|\leq \frac{\delta}{2}$ and $|z_2|\leq \frac{\delta}{2}$ from $T_1$ and $T_2$, respectively and glue the remaining parts along $\frac{\delta}{2}<|z_1|< \sqrt{\delta}$ and $\frac{\delta}{2}<|z_2|< \sqrt{\delta}$ by the formula $z_1.\overline{z_2}= \frac{\delta ^2}{2}$. Here the norms come from the induced metrics and to do gluing, we will assume that the induced metrics on $L_1$ and $L_2$ are compatible in the intervals $|z_1|< \sqrt{\delta}$ and $|z_2|< \sqrt{\delta}$. The approximate special Lagrangian submanifold will be locally a product manifold $H_\delta=T_{\delta}\times K$ where $T_{\delta}$ is the smoothing of $T_1\wedge T_2$ with respect to the gluing parameter $\delta$.

\vspace{.1in} 


More precisely, we can think of K as a curve in ${\bf C^3}$ with tangent ${\mathcal T}$. (${\mathcal T}_b$ being the value at
point b.) Let $\mathcal N_b$ be the complex two dimensional subspace of ${\bf C^3}$
orthogonal to ${\mathcal T}_b$ and $J.{\mathcal T}_b$. Then starting with some complex basis
$V_{1_{b_0}}, V_{2_{b_0}}$ for one of the vector spaces ${\mathcal N}_b$ (say ${\mathcal N}_{b_{0}}$), use orthogonal
projection to define nearby complex bases $V_{1_{b}}$, $V_{2_{b}}$ for ${\mathcal N}_b$. Then
one gets a real five dimensional manifold locally the image of

\vspace{.1in}

\hspace{.6in} ${\bf R}\times {\bf C} \times {\bf C} \rightarrow {\bf C^3}$ 

\vspace{.04in} 

\noindent by

\vspace{.04in} 

\hspace{.6in} $(t, z_1, z_2) \rightarrow c(t) + z_1 V_{1_{c(t)}} + z_2 V_{2_{c(t)}}$ 

\vspace{.1in} 

\noindent where $c(t) \in {\bf C^3}$ is the image point on the curve K under the
parameterization $c : {\bf R} \rightarrow K \subset {\bf C^3}$. Within this real 5 dimensional
manifold one has a real 3 dimensional manifold $H_\delta$ where $z_1 \overline {z_2} = \frac{\beta\delta^2}{2}$ for some cutoff function $\beta$.



\vspace{.1in}

When we do the smoothing for $T_1\wedge T_2$ the new surface can be written as a graph over $T_1$ and similarly over $T_2$. Therefore, locally, one can write the {\em approximate} special Lagrangian submanifold $H_\delta=L_1\# _\delta L_2(\kappa, z_1, z_2)$ for $\kappa\in K$ and $z_1\in \mathcal D_1$, $z_2\in \mathcal D_2$ :

\vspace{.2in}

$$L_1\# _\delta L_2(\kappa, z_1, z_2)= \left\{ \begin{array}{ccc}

(\kappa, 0, z_2) ,& \mbox{if} & |z_2| > \sqrt{\delta}\\

\\

(\kappa, \frac {\beta (z_2) {\delta }^2}{2\overline {z_2}}, z_2) ,& \mbox{if} & |z_1| < \frac{\delta}{2}$ and $|z_2| < \sqrt{\delta}\\

\\

(\kappa, \frac {{\delta }^2}{2\overline {z_2}}, z_2) ,& \mbox{if} &  \frac{\delta}{2}<|z_1| < \frac{\delta}{\sqrt 2}, \frac{\delta}{\sqrt 2}< |z_2| < \delta \\

\\

(\kappa, z_1, \frac {{\delta }^2}{2\overline {z_1}}) ,& \mbox{if} &  \frac{\delta}{\sqrt 2}< |z_1| < \delta, \frac{\delta}{2}< |z_2| < \frac{\delta}{\sqrt 2} \\

\\

(\kappa, z_1, \frac {\beta (z_1) {\delta }^2}{2\overline {z_1}}) ,& \mbox{if} &  \delta < |z_1|$ and $|z_2| < \frac{\delta}{2}\\

\\

(\kappa, z_1, 0) ,& \mbox{if} & \sqrt{\delta}<|z_1|\\

\end{array} \right.  $$

\vspace{.1in} 

\noindent where $\beta(z_i)$ is a smooth cut-off function supported in $[\frac{\delta}{\sqrt 2}, \sqrt{\delta}]$ with $\beta(z_i)$ being equal to 0 when $|z_i|> \sqrt{\delta} $, and $\beta(z_i)$ being equal to 1 when $|z_i|< \frac{\delta}{\sqrt 2} $. Moreover we will assume that $\beta(z_i)$ is radially symmetric with respect to $z_i$. This means that $\beta$ satisfies $x\frac{\partial\beta(z_1)}{\partial y}=y\frac{\partial\beta(z_1)}{\partial x}$ where the complex coordinate $z_1$ is represented by $(x,y)$ for $x,y \in {\bf R}$ and similarly $u\frac{\partial\beta(z_2)}{\partial v}=v\frac{\partial\beta(z_2)}{\partial u}$ for $z_2=(u,v)$, $u,v \in {\bf R}$. Note that the function $\beta(z_i) = \frac{log|z_i/\sqrt{\delta}|}{log\sqrt{\delta}}$ satisfies the necessary properties and after smoothing we can use this as the cutoff function.

\vspace{.1in}

Note that this construction is symmetric around $|z_1|=\frac{\delta}{\sqrt 2}$, $|z_2|=\frac{\delta}{\sqrt 2}$ and therefore it will be sufficient to investigate the other properties on one side. 

\vspace{.1in} 

This construction also shows that if we take $L_1$ as the graph of a function $f$ and $L_2$ as the graph of another function $g$ then $L_1\#_\delta L_2$ can be viewed as the graph of $\beta_\delta(z_2) f+\beta_\delta (z_1)g$ which depends on the choice of the gluing parameter $\delta$.

\vspace{.1in} 






\vspace{.1in}


The map $H_\delta$ is not special Lagrangian however we will see that the error term converges to zero as $\delta\rightarrow 0$ and our goal is to construct a {\em true} special Lagrangian submanifold near $H_\delta$. While the special Lagrangian equations are non-linear, it will be enough to find a solution of the corresponding linearized equation with estimates on the solutions which are independent of $\delta$. Here $\delta$ is the parameter defining the conformal structure of the connected sum and we will assume that it is small in our constructions.


\subsection {Lagrangian Property}

After constructing the approximate special Lagrangian submanifold the next step is to show that $H_\delta=L_1\#_\delta L_2$ is Lagrangian. In other words the symplectic 2-form $\omega$ on the ambient Calabi-Yau manifold $X^{2n}$ restricts to zero on $L_1\#_\delta L_2$. Therefore we need to introduce the symplectic structure in a small neighbourhood of the singular set.

\vspace{.1in}

Let $x_1=x, x_2=y, x_3=\kappa, x_4=\zeta, x_5=u, x_6=v$ be the local coordinates in an open ball $\tilde{U}$ around $b\in L\subset X^{6}$. Also let $x_1,x_2,x_3$ be the local coordinates for $L_1$ and $x_{4},x_5,x_{6}$ be the local coordinates for $L_2$. The basis $e_1,...,e_{6}$ will generate the the tangent space, ${\mathcal T}_bX$. Then we can take the dual basis $\omega_1,...\omega_{6}$ for the basis $e_1,...,e_{6}$ of ${\mathcal T}_bX^{6}$ and without loss of generality we can assume that the restriction of $e_1,e_2,e_{3}$ will be the basis for ${\mathcal T}_bL_1$ and similarly $e_{4},e_5,e_{6}$ will be the basis for ${\mathcal T}_bL_2$. Then we can find the dual basis $\omega_1,\omega_2,\omega_3$ for ${\mathcal T}^*_bL_1$ and $\omega_{4},\omega_{5},\omega_{6}$ for ${\mathcal T}^*_bL_2$ with respect to the induced metric.

\vspace{.1in}

Recall that we have the local model:

\vspace{.1in}

$L_1\cong K \times S^1\times S^1 \cong K \times T_1$

$L_2\cong K \times S^1\times S^1 \cong K \times T_2$.

\vspace{.1in}

Therefore, without loss of generality we can assume that $\omega_1,\omega_2\in {\mathcal T}^*(T_1)$, $\omega_{5},\omega_{6}\in {\mathcal T}^*(T_2)$ and  $\omega_3,\omega_4\in {\mathcal T}^*K $. Then in local coordinates, the symplectic structure in a small neighbourhood of the singularity is given as follows:

\vspace{.1in}

$\omega=\omega_1\wedge \omega_{5}+\omega_2\wedge \omega_{6}+\omega_3\wedge \omega_{4}$


{\lem : The connected sum $L_1 \#_\delta L_2$ is a Lagrangian submanifold of $X$ with respect to the symplectic structure defined as above.}

\vspace{.1in}

{\bf Proof:} We will show that the symplectic two form $\omega$ restricts to zero on $L_1 \#_\delta L_2$ for any $\delta$ (i.e $\omega|_{L_1 \#_\delta L_2}=0$). 

\vspace{.1in}

Recall that this construction is symmetric around $|z_1|=\frac{\delta}{\sqrt 2}$, $|z_2|=\frac{\delta}{\sqrt 2}$ and therefore it will be sufficient to prove the lemma for only one side. So here we will work with the interval $\frac{\delta}{\sqrt 2}<|z_1|$.

\vspace{.1in}

Suppose $L=$graph $f=\{x+if(x):x\in V\subset {\bf R^n}\}\subset {\bf R^n}\bigoplus i{\bf R^n}$ is the graph of a smooth function $f:V\rightarrow {\bf R^n}$, for some open set $V$ in ${\bf R^n}$. Let $\frac{\partial}{\partial x_i}$ be the standard basis for ${\bf R^n}$. Then the tangent space of $L$ at any point is generated by the following vectors:

\vspace{.1in}

$e_1+i\frac{\partial f}{\partial x_1},...,e_n+i\frac{\partial f}{\partial x_n}$, where $e_i=\frac{\partial}{\partial x_i}$.

\vspace{.1in}

Similarly the tangent space of the approximate special Lagrangian $H_{\delta}$ is generated by $3$ vectors.

\vspace{.1in}

$E_1= e_1+i(\frac{\partial}{\partial x}(\frac{{\delta}^2\beta(x,y)x}{2(x^2+y^2)})e_{5}+\frac{\partial}{\partial y}(\frac{{\delta}^2\beta(x,y)x}{2(x^2+y^2)}e_{6}))$

\vspace{.1in}

$E_2= e_2+i(\frac{\partial}{\partial x}(\frac{{\delta}^2\beta(x,y)y}{2(x^2+y^2)})e_{5}+\frac{\partial}{\partial y}(\frac{{\delta}^2\beta(x,y)y}{2(x^2+y^2)}e_{6}))$

\vspace{.1in}

$E_3=e_3$.

\vspace{.1in}

\noindent Then

\vspace{.1in}

$\omega|_{H_\delta}=(\omega_1\wedge \omega_{5}+\omega_2\wedge \omega_{6}+\omega_3\wedge \omega_{4})|_{H_\delta}$

\vspace{.1in}
$\omega (E_i,E_j)=0$ for all $i,j\geq 1,2$

\vspace{.1in}

$\omega (E_1,E_2)=(\frac{\partial}{\partial y}\beta(x,y))\frac{\delta^2x}{2(x^2+y^2)}-(\frac{\partial}{\partial x}\beta(x,y))\frac{\delta^2y}{2(x^2+y^2)}=0$. 

\vspace{.1in}

\noindent where we used the property of the cutoff function $x\frac{\partial\beta}{\partial y}=y\frac{\partial\beta}{\partial x}$. This calculation shows that the restriction of $\omega$ onto the tangent space of $H_\delta$ is zero for any $\delta$. This implies that $H_\delta$ is Lagrangian independent of $\delta$.

\vspace{.1in}
 
Recall that $L_1\#_\delta L_2$ can be viewed as the graph of $\beta_\delta(z_2) f+\beta_\delta (z_1)g$ where $L_1$ is the graph of a function $f$ and $L_2$ is the graph of another function $g$. Also we know that a Lagrangian submanifold of the cotangent bundle ${\mathcal T}^*W$ of a manifold $W$, is locally defined as the image of a section $d\Phi:W\rightarrow {\mathcal T}^*W$ for some function $\Phi: W\rightarrow {\bf R}$. Since $L_1$ and $L_2$ are both Lagrangian there are functions $F,G$ such that $f=\nabla F$ and $g=\nabla G$. Using the given cutoff function $\beta(z_i)$ and showing that $\beta=\nabla\tilde{\Phi}$ for some function $\tilde{\Phi}$, one can easily see that $L_1\#_\delta L_2$ is the graph of $\beta_\delta(z_2) f+\beta_\delta (z_1)g=\nabla\tilde{\Psi}$ for some function $\tilde{\Psi}$ which implies that $L_1\#_\delta L_2$ is also Lagrangian.

\subsection {Induced Metric and The Error Term}
For a positive integer $k$ and a number $1\leq p < \infty$ we can define the $L^p_{k,\delta}$-norm of a smooth function $F:\Omega \rightarrow {\bf R}$ by

\vspace{.1in}

\hspace{.1in}  $||F||_{L^p_{k,\delta}}=(\displaystyle\int_{\Omega}^{}\mathop{\Sigma}\limits_{|\tilde{\nu}|\leq k}|\partial ^{\tilde{\nu}} F|^p dvol_{\delta})^{1/p}$

\vspace{.1in}

\noindent where $\tilde{\nu} = (\tilde{\nu}_1,...,\tilde{\nu}_n)$ is a multi-index and $|\tilde{\nu}|= \tilde{\nu}_1+...+\tilde{\nu}_n$. 

\vspace{.1in}

Lemma 2.1 says that the approximate special Lagrangian $L_1\#_\delta L_2$ is Lagrangian and the error term for being a special Lagrangian will be $Im(\xi)|_{L_1\#_\delta  L_2}$. In this section we will show that $||$error term$||_ {L^2_{2,\delta}}
$ approaches to zero as the gluing parameter $\delta\rightarrow 0$. In our estimates we will always use $C$ to denote a uniform constant independent of the gluing parameter $\delta$ but its actual value may vary in different places. Also note that since there is symmetry in our construction, it is sufficient to show this for only the interval $|z_1|>\frac{\delta}{\sqrt 2}$ where the tangent space is generated by the vectors $E_1,E_2,E_3$ as in Lemma 2.1.

\vspace{.1in}

Since $(-1)^{n(n-1)/2}(i/2)^n\xi\wedge\overline{\xi}$ = $\omega^n/n!$ for any $n$, we can write the complex 3-form $\xi$ in local coordinates as follows:

\vspace{.1in}

$\xi=(\omega_1+i\omega_{5})\wedge(\omega_2+i\omega_{6})\wedge(\omega_3+i\omega_{4})$

\vspace{.1in}





\vspace{.1in}

Again as in Lemma 2.1, we will use the vectors $E_1,E_2,E_3$ which span the tangent space of the approximate special Lagrangian. Then we get

\vspace{.1in}








\vspace{.1in}





$Im(\xi)(E_1,E_2,E_3)= \frac{\partial}{\partial y}(\frac{\delta^2\beta(x,y)y}{2(x^2+y^2)})+\frac{\partial}{\partial x}(\frac{\delta^2\beta(x,y)x}{2(x^2+y^2)})$

\vspace{.1in}

$=\frac{\partial\beta}{\partial y}\frac{\delta^2y}{2(x^2+y^2)}+\beta\frac{\partial}{\partial y}(\frac{\delta^2y}{2(x^2+y^2)})+\frac{\partial\beta}{\partial x}\frac{\delta^2x}{2(x^2+y^2)}+\beta\frac{\partial}{\partial x}(\frac{\delta^2x}{2(x^2+y^2)})$

\vspace{.1in}

$=\frac{\partial\beta}{\partial y}\frac{\delta^2y}{2(x^2+y^2)}+\frac{\partial\beta}{\partial x}\frac{\delta^2x}{2(x^2+y^2)}$

\vspace{.1in}

Note that $\beta\frac{\partial}{\partial y}(\frac{\delta^2y}{2(x^2+y^2)})+
\beta\frac{\partial}{\partial x}(\frac{\delta^2x}{2(x^2+y^2)})=0$

\vspace{.1in}


Therefore the error term is given as $\frac{\partial\beta}{\partial y}\frac{\delta^2y}{2(x^2+y^2)}+\frac{\partial\beta}{\partial x}\frac{\delta^2x}{2(x^2+y^2)}$. This implies that for the intervals where
$\beta=0$ or $\beta=1$, the error term will be zero. In other words, the approximate special Lagrangian will be a true special Lagrangian for these intervals.

\vspace{.1in}

Next we will show that the value of the $L^2_{2,\delta}$ norm of the error term approaches to zero as $\delta\rightarrow 0$. Note that since the volume form depends on the metric the $L^2_{2,\delta}$ norm also depends on the induced metric.



\vspace{.1in} 
Recall that at $\mathcal Y$ the singularity is smoothed by the formula $z_1.\overline{z_2}= \frac{\delta ^2}{2}$ where $z_1=(x,y)$ and $z_2=(u,v)$ are the complex valued coordinates on $T_1$ and $T_2$, respectively and $K$ is the codimension two singularity. Since we work with Lagrangian submanifolds here we will switch back to real coordinates x,y and u,v. 

\vspace{.1in} 

 $z_1.\overline{z_2}= \frac{\delta ^2}{2}$ implies that $u= \frac {\delta^2 x}{2(x^2+y^2)}$ and $v=\frac {\delta^2 y}{2(x^2+y^2)}$. 

\vspace{.1in}

We will find the induced metric $\tilde{g}_{\delta}$ on $H_{\delta}=T_{\delta}\times K$ where $T_{\delta}$ is the smoothing of $T_1\wedge T_2$ with gluing parameter ${\delta}$. 

\vspace{.1in} 

$\tilde{g}_{\delta}=(dxdx+dydy+dudu+dvdv+d\kappa d\kappa+d\zeta d\zeta)|_{T_{\delta}\times K}$.

\vspace{.1in} 

For the intervals where the cutoff function $\beta=1$ we have

\vspace{.1in} 

$du=\frac {\delta^2(-x^2+y^2)}{2(x^2+y^2)^2}dx- \frac{\delta^2 2xy}{2(x^2+y^2)^2}dy$, and $dv=\frac {\delta^2(x^2-y^2)}{2(x^2+y^2)^2}dy- \frac{\delta^2 2xy}{ 2(x^2+y^2)^2}dx$  

\vspace{.1in} 

The induced metric for $\beta=1$ will be of the form

\vspace{.1in}


$\tilde{g}_{\delta}|_{H_{\delta}}=(1+\frac {\delta^4}{4(x^2+y^2)^2})dxdx+(1+\frac {\delta^4}{4(x^2+y^2)^2})dydy+ \mathcal Ad\kappa d\kappa$

\vspace{.1in}

\noindent where $\mathcal A$ is some small but not necessarily a constant function. If $\mathcal A=1$ this implies that our model is a geometric product. Recall that the $dvol_{\delta}$ is given by $\sqrt {det \tilde{g}_{\delta,ij}} dx_1...dx_n$ and we can easily calculate this term since we now know the induced metric and how it depends on $\delta$. 

\vspace{.1in}

$\sqrt {det \tilde{g}_{\delta,ij}}=(\mathcal A(1+\frac {\delta^4}{4(x^2+y^2)^2})(1+\frac {\delta^4}{4(x^2+y^2)^2}))^{1/2}=\sqrt{\mathcal A}(1+\frac {\delta^4}{4(x^2+y^2)^2})$

\vspace{.1in}

The next step is to find the volume form in general so that we can show that $L^2_{2,\delta}$ norm of the error term goes to 0 as the gluing parameter goes to 0. For the other intervals we cannot ignore the cutoff function, in fact the interval $0<\beta<1$ is the region where the error term is nonzero. 
 
\vspace{.1in} 

$z_1.\overline{z_2}= \frac{\beta\delta ^2}{2}$ implies that $u= \frac {\beta\delta^2 x}{2(x^2+y^2)}$ and $v=\frac {\beta\delta^2 y}{2(x^2+y^2)}$.
 
\vspace{.1in} 

As we did before, we can find $du$ and $dv$ terms:

\vspace{.1in}

$du=[ \frac {\delta^2(-x^2+y^2)}{2(x^2+y^2)^2}\beta(x,y)+\frac {\delta^2x}{2(x^2+y^2)}\frac{\partial \beta}{\partial x}]dx+ [\frac{-\delta^2 2xy}{2(x^2+y^2)^2}\beta(x,y)+ \frac{\delta^2x}{2(x^2+y^2)}\frac{\partial \beta}{\partial y}]dy$ and 

\vspace{.1in} 

$dv=[\frac {\delta^2(x^2-y^2)}{2(x^2+y^2)^2}\beta(x,y)+\frac {\delta^2y}{2(x^2+y^2)}\frac{\partial \beta}{\partial y}]dy+[ \frac{-\delta^2 2xy}{ 2(x^2+y^2)^2}\beta(x,y)+\frac {\delta^2y}{2(x^2+y^2)}\frac{\partial \beta}{\partial x}]dx$

\vspace{.1in} 

For simplicity, we will write $du=\tilde{A}dx+\tilde{B}dy$ and $dv=\tilde{C}dx+\tilde{D}dy$. Then as before
the induced metric will be of the form

\vspace{.1in}

$\tilde{g}_{\delta}|_{T_{\delta}\times K}=(1+\tilde{A}^2+\tilde{C}^2)dxdx+(1+\tilde{B}^2+\tilde{D}^2)dydy+\mathcal A d\kappa d\kappa$

\vspace{.1in}

\noindent and 

\vspace{.1in}

$\sqrt {det \tilde{g}_{\delta,ij}}=(\mathcal A(1+\tilde{A}^2+\tilde{C}^2)(1+\tilde{B}^2+\tilde{D}^2))^{1/2}$



\vspace{.1in}

\noindent $=\sqrt{\mathcal A}(1+(\frac{\delta^2(-x^2+y^2)}{2(x^2+y^2)^2}\beta +\frac {\delta^2x}{2(x^2+y^2)}\frac{\partial \beta}{\partial x})^2+(\frac{-\delta^2 2xy}{ 2(x^2+y^2)^2}\beta +\frac {\delta^2y}{2(x^2+y^2)}\frac{\partial \beta}{\partial x})^2)\cdot$

\hspace{.6in} $(1+(\frac{-\delta^2 2xy}{2(x^2+y^2)^2}\beta +\frac{\delta^2x}{2(x^2+y^2)}\frac{\partial \beta}{\partial y})^2+(\frac {\delta^2(x^2-y^2)}{2(x^2+y^2)^2}\beta +\frac {\delta^2y}{2(x^2+y^2)}\frac{\partial \beta}{\partial y})^2)$

\vspace{.1in}

\noindent $=\sqrt{\mathcal A}[1+\frac{\delta^4}{4(x^2+y^2)^2}(\beta^2-2x\beta\frac{\partial\beta}{\partial x}+(x^2+y^2)(\frac{\partial\beta}{\partial x})^2)]\cdot$

\hspace{.7in}$ [1+\frac{\delta^4}{4(x^2+y^2)^2}(\beta^2-2y\beta\frac{\partial\beta}{\partial y}+(x^2+y^2)(\frac{\partial\beta}{\partial y})^2)]$

\vspace{.1in}

\noindent $\leq\sqrt{\mathcal A}[1+\frac{\delta^4}{4(x^2+y^2)^2}((1-x\beta_x)^2+y^2\beta_x^2)]\cdot [1+\frac{\delta^4}{4(x^2+y^2)^2}((1-y\beta_y)^2+x^2\beta_y^2)]$

\vspace{.2in}

\noindent where $\beta_x=\frac{\partial\beta}{\partial x}$ and $\beta_y=\frac{\partial\beta}{\partial y}$. Using the basic properties of $\beta$, we can estimate the terms $(1-x\beta_x)^2+y^2\beta_x^2$ and $(1-y\beta_y)^2+x^2\beta_y^2$ by the term $(1-\frac{1}{\log\sqrt \delta})^2+\frac{1}{\log ^2 \sqrt \delta} $ and 
$(y\beta_y+x\beta_x)^2$ by $\frac{1}{\log ^2\sqrt \delta}$.



\vspace{.2in}

Then for a small ball $B_\delta$ around the singular point, $||$error term$||^2_{L^2_{\delta}}=||Im(\xi)(E_1,E_2,E_3)||^2_{L^2_{\delta}}$ can be bounded by the integral

\vspace{.1in}

\noindent ${\displaystyle\int_{B_\delta}}\sqrt{\mathcal A}\frac{\delta^4}{4(x^2+y^2)^2}\cdot\frac{1}{\log ^2\sqrt \delta}\cdot (1+\frac{\delta^4}{4(x^2+y^2)^2}(1-\frac{2}{\log\sqrt \delta}+\frac{2}{\log ^2 \sqrt \delta}))dxdyd{\kappa}$




  


\vspace{.1in}

\noindent $={\frac{1}{\log ^2\sqrt \delta}\displaystyle\int_{B_\delta}}\sqrt{\mathcal A}\frac{\delta^4}{4(x^2+y^2)^2}dxdyd\kappa +{\displaystyle\int_{B_\delta}}\sqrt{\mathcal A}\frac{\delta^8}{16(x^2+y^2)^4}\cdot(1-\frac{2}{\log\sqrt \delta}- \frac{2}{\log ^2\sqrt \delta})dxdyd\kappa$

\vspace{.1in}

\noindent $=\frac{\delta^4}{4\log ^2\sqrt \delta}[{\displaystyle\int_{B_\delta}}\sqrt{\mathcal A}\frac{1}{(x^2+y^2)^2}dxdyd\kappa+ {\displaystyle\int_{B_\delta}}\sqrt{\mathcal A}\frac{\delta^4}{4(x^2+y^2)^4}dxdyd\kappa$

\vspace{.1in}

\hspace{.5in}$-2{\displaystyle\int_{B_\delta}}\sqrt{\mathcal A}\frac{\delta^4}{4(x^2+y^2)^4}\frac{1}{\log\sqrt \delta}dxdyd\kappa-2{\displaystyle\int_{B_\delta}}\sqrt{\mathcal A}\frac{\delta^4}{4(x^2+y^2)^4}\frac{1}{\log ^2\sqrt \delta}dxdyd\kappa]$


\vspace{.1in}


\vspace{.1in}

\noindent where $B_{\delta}=\{\delta^2<x^2+y^2<\delta \}\times R\subset R^3$ and has coordinates $x,y,\kappa$. Switching to polar cordinates for $x,y$ and integrating each term separately over $B_{\delta}$ (for $0\leq\theta\leq 2\pi$ and $\delta\leq r\leq \sqrt{\delta})$ we get

\vspace{.1in}

\noindent $||$error term$||^2_{L^2_\delta}=|| (Im(\xi)(E_1,E_2,E_3))||^2_{L^2_{\delta}}$

\vspace{.1in}
 
$\leq C[\frac{\delta^3}{\log ^2\sqrt \delta}+\frac{\delta^2}{\log ^2\sqrt \delta}+\frac{\delta^5}{\log ^2\sqrt \delta}+\frac{\delta^5}{\log^3\sqrt \delta}+\frac{\delta^2}{\log ^3\sqrt \delta}+\frac{\delta^5}{\log ^4\sqrt \delta}+\frac{\delta^2}{\log ^4\sqrt \delta}]$

\vspace{.1in}

\noindent where $C$ is a constant independent of $\delta$. Here we also used the fact that $L$ is compact and $\displaystyle\int_{B_\delta}^{} \sqrt{\mathcal A} dxdyd\kappa$ can be bounded by some constant which does not depend on $\delta$.



\vspace{.1in}

The error term is a combination of terms which approach to $0$ with different rates as $\delta$ goes to 0. It is sufficient to bound the error term with $\frac{\delta^2}{\log ^2\sqrt \delta}$ which is the slowest term in this combination. So we can write $||$error term$||^2_{L^2_\delta}\leq C\frac{\delta^2}{\log ^2\sqrt \delta}\leq C\delta^2$. Next, we will estimate $||$error term$||_{L^2_{2,\delta}}$ by calculating $||\nabla (Im(\xi)(E_1,E_2,E_3))||^2_{L^2_{\delta}}$ and $||\nabla ^2(Im(\xi)(E_1,E_2,E_3))||^2_{L^2_{\delta}}$. Similar calculations and writing 

\vspace{.1in}

\noindent $\sqrt{\mathcal A}(1+\frac{\delta^4}{4(x^2+y^2)^2}(1-\frac{2}{\log\sqrt \delta}+\frac{2}{\log ^2 \sqrt \delta}))dxdyd{\kappa}$ for dvol$_\delta$ gives

\vspace{.1in}

\noindent $||\nabla $(error term)$||^2_{L^2_{\delta}}=\frac{\delta^4}{4} {\displaystyle\int_{B_\delta}} |\nabla (\frac{y\beta_y+x\beta_x}{x^2+y^2})|^2 $dvol$_\delta $ 



\vspace{.1in}

\hspace{.15in} $\leq \frac{\delta^4}{\log ^2 \sqrt \delta} {\displaystyle\int_{B_\delta}}\sqrt{\mathcal A}(\frac{1}{x^2+y^2})\cdot(1+\frac{\delta^4}{4(x^2+y^2)^2}(1-\frac{2}{\log\sqrt \delta}+\frac{2}{\log ^2 \sqrt \delta}))dxdyd{\kappa}$ 


\vspace{.1in}

\hspace{.15in} $\leq C[\frac{\delta^4}{\log \sqrt \delta}+\frac{\delta^6}{\log ^2\sqrt \delta}+\frac{\delta^4}{\log ^2\sqrt \delta}+\frac{\delta^6}{\log^3\sqrt \delta}+\frac{\delta^4}{\log ^3\sqrt \delta}+\frac{\delta^6}{\log ^4\sqrt \delta}+\frac{\delta^4}{\log ^4\sqrt \delta}]$

\vspace{.1in}

\noindent and

\vspace{.1in}

\noindent $||\nabla ^2$(error term)$||^2_{L^2_{\delta}}=\frac{\delta^4}{4} {\displaystyle\int_{B_\delta}} |\nabla ^2 (\frac{y\beta_y+x\beta_x}{x^2+y^2})|^2 $dvol$_\delta $


\vspace{.1in}

\hspace{.15in} $\leq \frac{\delta^4}{\log ^2 \sqrt \delta} {\displaystyle\int_{B_\delta}}\sqrt{\mathcal A}(\frac{1}{(x^2+y^2)^2})\cdot(1+\frac{\delta^4}{4(x^2+y^2)^2}(1-\frac{2}{\log\sqrt \delta}+\frac{2}{\log ^2 \sqrt \delta}))dxdyd{\kappa}$ 


\vspace{.1in}

\hspace{.15in} $\leq C[\frac{\delta^3}{\log ^2\sqrt \delta}+\frac{\delta^2}{\log ^2\sqrt \delta}+\frac{\delta^5}{\log ^2\sqrt \delta}+\frac{\delta^5}{\log^3\sqrt \delta}+\frac{\delta^2}{\log ^3\sqrt \delta}+\frac{\delta^5}{\log ^4\sqrt \delta}+\frac{\delta^2}{\log ^4\sqrt \delta}]$

\vspace{.1in}

Combining all these estimates we can finally conclude the following lemma.


{\lem : There exists a constant $C$, independent of $\delta$ such that the following holds:

\vspace{.1in}

$||Im(\xi)(E_1,E_2,E_3)||_{L^2_{2,\delta}}\leq C.{\delta}$  and this error term converges to zero as $\delta\rightarrow 0$.}

\section{Special Lagrangian Equation}

In this section, we will study the submanifolds close to the approximate special Lagrangian submanifold $H_\delta=L_1\#_\delta L_2$. More precisely we will write a nonlinear partial differential equation such that the solution set gives the special Lagrangian submanifolds near $H_\delta$. Then we will study the linearized operator $D_{\delta}=D_{H_\delta }$ for this equation and show that it is Fredholm.

\vspace{.1in}
  
Since $H_\delta$ is compact there is a tubular neighbourhood of $H_\delta$ in $X$ which is identified via the normal exponential map to a neighbourhood $N_\gamma(H_\delta)=\{\mathcal V\in N(H_\delta)|||\mathcal V||<\gamma\}$ of the zero section in the normal bundle $N(H_\delta)$. Therefore we can identify nearby submanifolds with small vector fields.

\vspace{.1in}

Also for a small normal vector field $\mathcal V$ in $\Gamma(N(H_\delta))$, the space of sections of the normal bundle, we can define the deformation map for the approximate special Lagrangian $H_{\delta}$ as follows,

\vspace{.1in}     

${\mathcal F}: \Gamma(N(H_\delta))\rightarrow \Omega^2(H_\delta)\bigoplus\Omega^3(H_\delta)$

\vspace{.1in}

${\mathcal F}(\mathcal V)=((\exp_{\mathcal V})^*(-\omega), (\exp_{\mathcal V})^*(Im(\xi))$

\vspace{.1in}

The deformation map ${\mathcal F}$ is the restriction of $-\omega$ and $Im(\xi))$ to $(H_{\delta})_{\mathcal V}$ and then pulled back to $H_\delta$ via $(\exp_{\mathcal V})^*$. Then ${\mathcal F}^{-1}(0,0)$ will correspond to the space of special Lagrangian submanifolds near $H_\delta$. 

\vspace{.1in}

We showed in 2.2 that $H_\delta$ is a Lagrangian submanifold of $X$ independent of $\delta$. In order to simplify the calculations we will assume that the deformation vector field $\mathcal V$ preserves this property. It is well known that for a special Lagrangian submanifold the linearization of ${\mathcal F}$ at 0 is $d+*d^*$, [13]. Here we will show that the linearized operator of ${\mathcal F}$ at 0 for an approximate special Lagrangian is given by $d+*(\Psi \cdot d^*)$ where $\Psi$ is a small function which is equal to 1 in the intervals with zero error term. Since $H_\delta$ is deformed as Lagrangian we obtain $d$ as the first part of the linearization as before and therefore in this section we will only study the second part of ${\mathcal F}$ which depends on $\delta$.

\vspace{.1in}  

Again let $z_1\in T_1$ and $z_2\in T_2$ be the complex valued coordinates. Away from the singularity, for $\sqrt{\delta}\leq |z_1|$ and $\sqrt{\delta}\leq |z_2|$, $H_\delta$ can be identified with $L_1$ and $L_2$, respectively and therefore $H_\delta$ is a real special Lagrangian and the linearized operator is given by $d+*d^*$ as before. Recall that the error term is zero in the middle of the neck area where $\frac{\delta}{2}<|z_1| < \frac{\delta}{\sqrt 2}, \frac{\delta}{\sqrt 2}< |z_2| < \delta$ and also $\frac{\delta}{\sqrt 2}< |z_1| < \delta, \frac{\delta}{2}< |z_2| < \frac{\delta}{\sqrt 2}$. So it remains to investigate the special Lagrangian equation and in particular the linearized operator for the intervals $|z_1| < \frac{\delta}{2}, |z_2| < \sqrt{\delta}$ and $|z_2| < \frac{\delta}{2}, \delta < |z_1|$ where the error term is small but nonzero.









\vspace{.1in}  

First we will find a condition on the deformation vector field $\mathcal V$ such that the approximate special Lagrangian can be deformed as Lagrangian.

\vspace{.1in}  

A Lagrangian submanifold of the cotangent bundle $T^* W$ of a manifold $W$, (for which the projection to $W$ is a local diffeomorphism) is locally defined as the image of a section $d\Phi:W\rightarrow T^*W$ for some function ${\mathcal \Phi}: W\rightarrow {\bf R}$. This implies that the deformation vector field $\mathcal V$ should be of the form $\mathcal V=J$. grad$(\Phi)$ for some function $\Phi$ on the approximate special Lagrangian submanifold and an almost complex structure $J$. One can easily show that this is equivalent to saying that the one-form that corresponds to the deformation is an exact one form, i.e. $\mathcal V$ is a Hamiltonian deformation vector field.

\vspace{.1in}

One other advantage of assuming the deformation vector field $\mathcal V$ preserves Lagrangian property is that we will be able to switch from one forms to functions when we work with the second part of the linearized operator.

\vspace{.1in}

Next we will find the condition which gives nearby special Lagrangian submanifolds. Recall that we can view a nearby special Lagrangian submanifold as a graph of a smooth function f. The graph of $f$ is special Lagrangian if it is Lagrangian (i.e. $f=\nabla F$ for some scalar function $F:H_\delta\rightarrow {\bf R}$ and $Im\xi(e_1+i\frac{\partial f}{\partial x_1},...,e_3+i\frac{\partial f}{\partial x_3})$ is zero.



\vspace{.1in}

Recall that $z_1.\overline{z_2}= \frac{\beta\delta ^2}{2}$ implies that $u= \frac {\beta\delta^2 x}{2(x^2+y^2)}$ and $v=\frac {\beta\delta^2 y}{2(x^2+y^2)}$. Therefore on the approximate special Lagrangian submanifold we will use the coordinates $u,v$ instead of $x,y$. The tangent space of the approximated special Lagrangian is generated by $3$ vectors.

\vspace{.1in}

$E_1= e_1+i(\frac{\partial}{\partial x}(\frac{{\delta}^2\beta(x,y) x}{2(x^2+y^2)})e_{5}+\frac{\partial}{\partial x}(\frac{{\delta}^2\beta(x,y) y}{2(x^2+y^2)}e_{6}))$

\vspace{.1in}

$E_2= e_2+i(\frac{\partial}{\partial y}(\frac{{\delta}^2\beta(x,y)x}{2(x^2+y^2)})e_{5}+\frac{\partial}{\partial y}(\frac{{\delta}^2\beta(x,y)y}{2(x^2+y^2)}e_{6}))$

\vspace{.1in}

$E_3=e_3$

\vspace{.1in}

The set of vectors that span the tangent space of the nearby submanifolds which can be written as a graph of  a function $f$. Note that $f$ depends on the gluing parameter $\delta$ but for simplicity we will just write $f$ instead of $f_\delta$. They are given as  

\vspace{.1in}
$E'_1= E_1+i(\frac{\partial f_1}{\partial u} E_{5}+\frac{\partial f_1}{\partial v}E_{6})$, $E'_2= E_2+i(\frac{\partial f_2}{\partial u} E_{5}+\frac{\partial f_2}{\partial v}E_{6})$

\vspace{.1in}

$E'_3= E_3$

\vspace{.1in}

$E_{5}=JE_1=-\frac{\partial}{\partial x}(\frac{{\delta}^2\beta x}{2(x^2+y^2)})e_{1}-\frac{\partial}{\partial x}(\frac{{\delta}^2\beta y}{2(x^2+y^2)})e_{2}+e_{5}$

\vspace{.1in}

$E_{6}=JE_2=-\frac{\partial}{\partial y}(\frac{{\delta}^2\beta x}{2(x^2+y^2)})e_{1}-\frac{\partial}{\partial y}(\frac{{\delta}^2\beta y}{2(x^2+y^2)})e_{2}+e_{6}$

\vspace{.1in}

\noindent where $J$ is the almost complex structure. We restrict $Im\xi$ to this tangent space 

\vspace{.1in}

\noindent $E'_1=e_1+\frac{\partial}{\partial x}(\frac{{\delta}^2\beta x}{2(x^2+y^2)})e_{5}+\frac{\partial}{\partial x}(\frac{{\delta}^2\beta y}{2(x^2+y^2)})e_{6}+\frac{\partial f_1}{\partial u}(-\frac{\partial}{\partial x}(\frac{{\delta}^2\beta x}{2(x^2+y^2)})e_{1}$

\hspace{.05in} $-\frac{\partial}{\partial x}(\frac{{\delta}^2\beta y}{2(x^2+y^2)})e_{2}+e_{5})+\frac{\partial f_2}{\partial u}(-\frac{\partial}{\partial y}(\frac{{\delta}^2\beta x}{2(x^2+y^2)})e_{1}-\frac{\partial}{\partial y}(\frac{{\delta}^2\beta y}{2(x^2+y^2)})e_{2}+e_{6})$

\vspace{.1in}

\vspace{.1in}

\noindent $E'_2=e_2+\frac{\partial}{\partial y}(\frac{{\delta}^2\beta x}{2(x^2+y^2)})e_{5}+\frac{\partial}{\partial y}(\frac{{\delta}^2\beta y}{2(x^2+y^2)})e_{6}+\frac{\partial f_1}{\partial v}(-\frac{\partial}{\partial x}(\frac{{\delta}^2\beta x}{2(x^2+y^2)})e_{1}$

\hspace{.05in} $-\frac{\partial}{\partial x}(\frac{{\delta}^2\beta y}{2(x^2+y^2)})e_{2}+e_{5})+\frac{\partial f_2}{\partial v}(-\frac{\partial}{\partial y}(\frac{{\delta}^2\beta x}{2(x^2+y^2)})e_{1}-\frac{\partial}{\partial y}(\frac{{\delta}^2\beta y}{2(x^2+y^2)})e_{2}+e_{6})$

\vspace{.1in}

\noindent $E'_3= E_3$

\vspace{.1in}

\noindent Then we get 

\vspace{.1in}

\noindent $Im\xi(E'_1,E'_2,E'_3)=Im[(\omega_1+i\omega_{5})\wedge(\omega_2+i\omega_{6})\wedge(\omega_3+i\omega_{4})](E'_1,E'_2,E'_3)$

\vspace{.1in}
\noindent =$[1+(\frac{\partial f_1}{\partial u})(\frac{\partial}{\partial x}(-\frac{{\delta}^2\beta x}{2(x^2+y^2)}))+(\frac{\partial f_2}{\partial u})(\frac{\partial}{\partial y}(-\frac{{\delta}^2\beta x}{2(x^2+y^2)}))]\cdot[(\frac{\partial}{\partial y}(\frac{{\delta}^2\beta y}{2(x^2+y^2)}))+(\frac{\partial f_2}{\partial v})]$
\vspace{.1in}

$-[\frac{\partial}{\partial x}(\frac{{\delta}^2\beta y}{2(x^2+y^2)})+(\frac{\partial f_2}{\partial u})]\cdot[(\frac{\partial f_1}{\partial v})(\frac{\partial}{\partial x}(-\frac{{\delta}^2\beta x}{2(x^2+y^2)}))+(\frac{\partial f_2}{\partial v})(\frac{\partial}{\partial y}(-\frac{{\delta}^2\beta x}{2(x^2+y^2)}))]$

\vspace{.1in}
$+[\frac{\partial}{\partial x}(\frac{{\delta}^2\beta x}{2(x^2+y^2)})+(\frac{\partial f_1}{\partial u})]\cdot[1+(\frac{\partial f_1}{\partial v})(\frac{\partial}{\partial x}(-\frac{{\delta}^2\beta y}{2(x^2+y^2)}))+(\frac{\partial f_2}{\partial v})(\frac{\partial}{\partial y}(-\frac{{\delta}^2\beta y}{2(x^2+y^2)}))]$

\vspace{.1in}

$ -[\frac{\partial}{\partial y}(\frac{{\delta}^2\beta x}{2(x^2+y^2)})+(\frac{\partial f_1}{\partial v})]\cdot[(\frac{\partial f_1}{\partial u})(\frac{\partial}{\partial x}(-\frac{{\delta}^2\beta y}{2(x^2+y^2)}))+(\frac{\partial f_2}{\partial u})(\frac{\partial}{\partial y}(-\frac{{\delta}^2\beta y}{2(x^2+y^2)}))]$

\vspace{.1in}

\noindent $=\frac{\partial}{\partial x}(\frac{{\delta}^2\beta x}{2(x^2+y^2)})+\frac{\partial}{\partial y}(\frac{{\delta}^2\beta y}{2(x^2+y^2)})+(\frac{\partial f_1}{\partial u})+(\frac{\partial f_2}{\partial v})$

\vspace{.1in}

$+[(\frac{\partial f_1}{\partial u})+(\frac{\partial f_2}{\partial v})](-(\frac{\partial}{\partial x}(\frac{{\delta}^2\beta x}{2(x^2+y^2)}))(\frac{\partial}{\partial y}(\frac{{\delta}^2\beta y}{2(x^2+y^2)}))+(\frac{\partial}{\partial y}(\frac{{\delta}^2\beta x}{2(x^2+y^2)}))(\frac{\partial}{\partial x}(\frac{{\delta}^2\beta y}{2(x^2+y^2)}))$

\vspace{.1in}

$-[(\frac{\partial f_1}{\partial u})(\frac{\partial f_2}{\partial v})-(\frac{\partial f_1}{\partial v})(\frac{\partial f_2}{\partial u})](\frac{\partial}{\partial x}(\frac{{\delta}^2\beta x}{2(x^2+y^2)})+\frac{\partial}{\partial y}(\frac{{\delta}^2\beta y}{2(x^2+y^2)}))$

\vspace{.1in}

The special Lagrangian equation is

\vspace{.1in}

$\Delta \mathcal G +(1-$detHess $\mathcal G)\Delta F+ \Delta \mathcal G(\frac{\partial^2F}{\partial u \partial v}
\frac{\partial^2F}{\partial v \partial u} - \frac{\partial^2F}{\partial u \partial u}\frac{\partial^2F}{\partial v \partial v})=0$.

\vspace{.1in}

\noindent where $\nabla F_\delta=f_\delta$ and $\frac{\partial \mathcal G}{\partial x}= \frac{{\delta}^2\beta x}{2(x^2+y^2)},\frac{\partial \mathcal G}{\partial y}= \frac{{\delta}^2\beta y}{2(x^2+y^2)}$. 

\vspace{.1in}

Note that in this equation the operator $\Delta \mathcal G$ is given in terms of $x,y$ variables and  $\Delta F$ is given in terms of $u,v$. The linearization gives $(1-$detHess $\mathcal G)\Delta F$ which is equivalent to $\Delta F$ for the intervals where the cutoff function $\beta =0$ and $\beta =1$. 

\vspace{.1in}

One can show that ${\displaystyle\int_{B_\delta}}|$detHess$ \mathcal G_\delta|^2 $dvol$_{\delta}$

\vspace{.1in}

\noindent $={\displaystyle\int_{B_\delta}}|(-(\frac{\partial}{\partial x}(\frac{{\delta}^2\beta x}{2(x^2+y^2)}))(\frac{\partial}{\partial y}(\frac{{\delta}^2\beta y}{2(x^2+y^2)}))+(\frac{\partial}{\partial y}(\frac{{\delta}^2\beta x}{2(x^2+y^2)}))(\frac{\partial}{\partial x}(\frac{{\delta}^2\beta y}{2(x^2+y^2)}))|^2$dvol$_\delta$

\vspace{.1in}

\noindent $\leq \frac{\delta ^8}{16}{\displaystyle\int_{B_\delta}}|\frac{1}{(x^2+y^2)^2}-\frac{1}{\log \sqrt \delta (x^2+y^2)^2}|^2$dvol$_\delta$

\vspace{.1in}

\noindent $\leq \frac{\delta ^8}{16}{\displaystyle\int_{B_\delta}}|\frac{1}{(x^2+y^2)^2}-\frac{1}{\log \sqrt \delta (x^2+y^2)^2}|^2 \cdot\sqrt{\mathcal A}(1+\frac{\delta^4}{4(x^2+y^2)^2}(1-\frac{2}{\log\sqrt \delta}+\frac{2}{\log ^2 \sqrt \delta}))dxdyd{\kappa}$

\vspace{.1in}

\noindent $\leq C\delta ^2 $

\vspace{.1in}

Therefore ${\displaystyle\int_{B_\delta}}|$detHess $ \mathcal G_\delta|^2 $dvol$_{\delta}$ is uniformly bounded and for small values of $\delta$, $||$detHess$ \mathcal G_\delta||$ will be small. For $0<\beta<1$ we have




\vspace{.1in}

\noindent $||((1-$detHess$ \mathcal G_{\delta})(\Delta_{\delta} ))^{-1}||=||(\Delta^{-1}_{\delta})(1-$detHess$ \mathcal G_\delta)^{-1}||$

\vspace{.1in}

\hspace{1.75in} $\leq ||\Delta_{\delta}^{-1}||\cdot||(1-$detHess$ \mathcal G_\delta)^{-1}||$  

\vspace{.1in}

\noindent and for sufficiently small $\delta$, we can take $||$detHess$ \mathcal G _{\delta}||$ to be arbitrarily small. Hence to invert the operator $(1-$detHess$ \mathcal G_{\delta})\Delta_{\delta}$ it is sufficient to check the invertibility of $\Delta_{\delta} $. In the next section, we will show this by the eigenvalue estimates.


\section{An Eigenvalue Estimate for the Linearized Operator}

\vspace{.1in}

In section 3 we showed that the second part of the linearized operator is given by $\Psi \cdot (d^*d)$ which depends on $\delta$ and to invert this it is sufficient to check the invertibility of $\Delta_{\delta} $. Next we obtain uniform estimates in the first eigenvalue of $\Delta_{\delta} $. More precisely we will show that for sufficiently small gluing parameter $\delta$, the linear operator $D_{\delta}=\Delta_{\delta} $ has a right inverse $Q_{\delta}$ which satisfies the following uniform inequality:

\vspace{.1in}

\hspace{.6in}  $||Q_{\delta} g||_{\delta,k+2,p}\leq C||g||_{\delta,k,p}$

\vspace{.1in}

\noindent for any $g \in L^p(H_\delta)$ and $k\leq 2$ where $C$ is independent of the gluing parameter $\delta$.   

\vspace{.1in}

The main estimate in this section is 

\vspace{.1in}

{\lem : There is a constant $C>0$ independent of the gluing parameter $\delta$, such that for $\delta$ sufficiently small, the first (nonzero) eigenvalue $\lambda_1(\Delta_{\delta})$ of $\Delta_{\delta}$ is bounded below by $C$.}

\vspace{.1in}

{\bf Proof:} We prove it by contradiction. Suppose that the lemma is not true. Then we may assume that the first eigenvalue $\lambda_1(\triangle_{\delta})$ converges to zero as $\delta$ tends to zero. Let $\phi_\delta$ be the eigenfunction of $\lambda_1(\Delta_{\delta})$ satisfying

\vspace{.1in}

$\displaystyle \int_{H_\delta} |\phi_\delta|^2 =1$ and $\displaystyle \int_{H_\delta} \phi_\delta =0$ and $\Delta_{\delta}\phi_\delta=\lambda_{1,\delta}\phi_\delta$ .

\vspace{.1in}

\noindent Here $\lambda_{1,\delta}$ determines the dependence of the first eigenvalue on the gluing parameter $\delta$ and ${H_\delta}$ is the smoothed approximate special Lagrangian.  Note that every time we change $\delta$, we change the induced metric on the approximated special Lagrangian and since the Laplacian operator depends on the metric, the eigenvalues of $\Delta_{\delta}$ depend on $\delta$. For simplicity we will use $\lambda_{\delta}$ for $\lambda_{1,\delta}$.

\vspace{.1in}

For small compact sets away from singularity, the $L^2_\delta$ norm is uniformly equivalent to the usual $L^2$ norm. On these compacts sets there exists a subsequence of $\phi_n$ that converges smoothly to a limit $\Delta \phi_0=0$. Following the same argument for the sequence of compact sets, and passing to a diagonal subsequence, we obtain a nonzero eigenfunction $\phi_0$ as the limit defined in the complement of the singularity satisfying 

\vspace{.1in}

$\displaystyle \int |\phi_0|^2 =1$ and $\displaystyle \int \phi_0 =0$

\vspace{.1in}

We now explain why $\phi_0$ cannot be zero.  If $\phi_0=0$ then for very small $\delta$, $\phi_{\delta}$ will be very small everywhere (almost zero) which contradicts the fact that 

\vspace{.1in}
$||\phi_{\delta}||_{L^2}\leq ||\phi_{\delta}||_{L^\infty}$ and our initial assumption $||\phi_{\delta}||_{L^2}=1$.

\vspace{.1in}

So we have a nonzero function $\phi_0$ in the limit and since $\lambda_{\delta}\rightarrow 0$ we get $\Delta_0\phi_0=0$.

\vspace{.1in}

On a compact manifold the only harmonic functions are constant functions. Therefore $\phi_0$ should be some nonzero constant. On one component $\phi_{\delta}$ will converge to a constant and on the other component it will converge to another constant. Since we assume that the singularity is irreducible these two constants should be same and since $\displaystyle \int \phi_0 =0$ this is only possible if $\phi_{\delta}$ converges to zero. This contradicts the fact that $\phi_0$ is nonzero. 

\vspace{.1in}

Next we should check the possibility where eigenfunctions are concentrating in the neck area. In this case the eigenfunctions can get trapped in the neck region and as the gluing parameter $\delta$ goes to 0 they converge to maps which are identically zero everywhere but blow up at one point. Recall that the neck region, the region close to the singularity is locally a product $T^2_{\delta}\times K\subset {\bf C^2}\times {\bf C}$ where $T_{\delta}$ is the neck area of $T_1\#_{\delta}T_2$ and $K$ is the codimension-2 singularity. Using compactly supported functions and Theorem 4.2 which is a well known fact from {\em Geometric Measure Theory} [17], we will show that we can also find a uniform lower bound for the first eigenvalue close to the singularity.




\vspace{.1in}

{\thm:  Suppose $h\in C_0^1(U^{n+k})$ for some domain $U$, and $h\geq0$.
 
\noindent Then for $p=\frac{n}{n-1}$, 

\vspace{.1in}

\hspace{1in} $(\displaystyle \int_N{h^p})^{1/p}\leq c (\int_N{|\nabla h|+h|{\mathcal H}|)}$

\vspace{.1in}

\noindent where $N$ is an n-dimensional domain and ${\mathcal H}$ is the mean curvature of N .}

\vspace{.1in}

In order to use Theorem 4.2 in our eigenvalue estimates we need to show that for any $N_{\delta}$, the $L^2$ norm of the mean curvature ${\mathcal H}_{\delta}$ is uniformly bounded and can be made arbitrarily small.

{\lem: For each small constant $C>0$ there is a $\delta$ so that the mean
curvature of $ N_{\delta}$ satisfies $||{\mathcal H}_\delta ||_{L^2} \leq C$}.

\vspace{.1in}


{\bf Proof:} It is easy to see that ${\mathcal H_\delta}$ is uniformly bounded. The derivative of the   
cutoff function $\beta$ is non-zero only in a region that is uniformly bounded away
from the singularity and in that region the approximate special Lagrangians
are converging in $C^1$ to the limit Lagrangian, which has ${\mathcal H}=0$.
Therefore the mean curvatures of the approximate special Lagrangians are converging to 0 uniformly. 


\vspace{.1in}

Next, we need to measure the deviation of the metric from the product metric on $T_\delta \times K$ and the contribution of the mean curvature of $K$ itself. Recall that locally we can work in ${\bf C^3}$ and since $K$ is one dimensional we can take $K=c(t)$ to be a curve in ${\bf C^3}$. Then one gets a real five dimensional manifold locally the image of



\vspace{.1in} 

\hspace{.6in} $\tilde{\theta}: {\bf R}\times {\bf C} \times {\bf C} \rightarrow {\bf C^3}$ 

\vspace{.04in}
\noindent by
\vspace{.04in} 

\hspace{.6in} $\tilde{\theta}(t, z_1, z_2) \rightarrow c(t) + z_1 V_1(t) + z_2 V_2(t)$ 

\vspace{.1in} 

\noindent where $V_1=\tilde{P}+i\tilde{R}$ and $V_2=\tilde{Q}+i\tilde{S}$ are complex valued normal vectors that span the complement of $c(t)$ in ${\bf C^3}$ and $z_1\in T_1$, ${z_2}\in T_2$. 


\vspace{.1in} 


Within this real 5 dimensional manifold one has a real 3 dimensional manifold which locally represents the approximate special Lagrangian as the image of $\tilde{l}$
where $\tilde{l}(t,z_1)=\tilde{\theta}(t, z_1, \frac{\delta^2}{\overline z_2})$. Note that $V_1$ and $V_2$ are changing as we move along $K$.

\vspace{.1in} 

Using this local model, one can calculate the mean curvature of $N_\delta$ as

\vspace{.1in} 

$\mathcal H_\delta = \displaystyle \frac{c''(t)+z_1V_1''(t)}{1+|z_1|^2}+ \delta^2(-\tilde{Z}c''(t)+\frac{V_2''(t)}{z_1(1+|z_1|^2)})+O(\delta^4)$

\vspace{.1in} 

\noindent where $\tilde{Z}$ is some bounded function. Here we will also assume that $|c''(t)|$, $|V_1''(t)|$ and $|V_1''(t)|$ are all bounded and do not depend on $\delta$. Then one can show that

\vspace{.1in}

$\displaystyle \int_{B_\delta} |\frac{c''(t)+z_1V_1''(t)}{1+|z_1|^2}|^2 $dvol$_\delta\leq C(\delta)$ where $C(\delta)\rightarrow 0$ as $\delta\rightarrow 0$.
 
\vspace{.1in}

Similar calculations for the other terms imply that $||\mathcal H_\delta||_{L^2}$ can be made arbitrarily small in the region $B_\delta=\{\delta^2<x^2+y^2<\delta\}\times R$. This is also true for $||\mathcal H_\delta||_{L^3}$. Showing that all Christoffel symbols are uniformly bounded will generalize this fact to other Calabi-Yau manifolds. This gives the proof of Lemma 4.3.

\vspace{.1in}






Next, we will use Theorem 4.2 and Lemma 4.3 to obtain a uniform lower bound for the first eigenvalue in $N_{\delta}$.


{\lem : For any $h\in C_0^1(N_{\delta})$ there exists a constant $C$, independent of $\delta$ such that the following inequality holds

\vspace{.1in}

\hspace{1.1 in} $C\leq \frac{||\nabla h||_{L^2}}{||h||_{L^2}}$  .}

\vspace{.1in}

{\bf Proof:} Recall from Theorem 4.2 that for $h\in C_0^1(N_{\delta})$ ,
we have the following inequality :

\vspace{.1in}

$(\displaystyle \int_{N_{\delta}}{h^{\frac{n}{n-1}})}^{\frac{n-1}{n}}\leq c (\displaystyle \int_{N_{\delta}}|\nabla h|+h|{\mathcal H}_{\delta}|)$

\vspace{.1in}

First we will show that for $n=2$ this inequality will imply $C\leq \frac{||\nabla h||_{L^2}}{||h||_{L^2}}$ and then generalize this to dimensions $n>2$.

\vspace{.1in}

$(\displaystyle \int_{N_{\delta}}{h^p})^{1/p}\leq c (\displaystyle \int_{N_{\delta}}{|\nabla h|+h|{\mathcal H}_{\delta}|)}=c (\displaystyle \int_{N_{\delta}}{|\nabla h|.1+h|{\mathcal H}_{\delta}|)}$

\vspace{.1in}
$\leq c (\displaystyle \int_{N_{\delta}}{|\nabla h|^2})^{1/2}(\displaystyle \int_{N_{\delta}}{1})^{1/2} +c(\displaystyle \int_{N_{\delta}}{|h|^2})^{1/2}(\displaystyle \int_{N_{\delta}}{|{\mathcal H}|^2})^{1/2}$

\vspace{.1in}

From Lemma 4.3 we know that the mean curvature is uniformly bounded and can be made small. So we can replace $||{\mathcal H}_{\delta}||_{L^2}$ with some small number $C_1$



\vspace{.1in}

$||h||_{L^2}\leq c (\displaystyle \int_{N_{\delta}}{|\nabla h|^2})^{1/2}(\displaystyle \int_{N_{\delta}}{1})^{1/2} + cC_1 (\displaystyle \int_{N_{\delta}}{|h|^2})^{1/2}$

\vspace{.1in}

\hspace{.42in} $\leq C ||\nabla h||_{L^2} +cC_1 ||h||_{L^2}$

\vspace{.1in}

Then $||h||_{L^2}\leq C ||\nabla h||_{L^2} +cC_1 ||h||_{L^2}$ and we get $\frac{1-cC_1}{C}\leq \frac{||\nabla h||_{L^2}}{||h||_{L^2}}$. This implies that for n=2 we get a uniform lower bound for the first eigenvalue in the neck area $N_{\delta}$. 

{\rem:} In Theorem 4.2, $n$ is the dimension of $N_\delta$. The theorem holds for p=2 only
when $\frac{n}{n-1} \geq 2$,  so only when n=1 or n=2. Since we work with a 3-dimensional Lagrangian we need to generalize this fact to higher dimensions.

\vspace{.1in}

Next we will generalize to $n>2$. First for any $h\in C_0^1(N_{\delta})$ we will replace $h$ with $h^\alpha$ for some $\alpha\in {\bf Q}$. Then $h^{\alpha}\in C_0^1(N_{\delta})$  and we have the following inequality :

\vspace{.1in}

$(\displaystyle \int_{N_{\delta}}{h^{\frac{\alpha n}{n-1}})}^{\frac{n-1}{n}}\leq c (\displaystyle \int_{N_{\delta}}|\nabla h^{\alpha}|+h^{\alpha}|{\mathcal H}_{\delta}|)$

\vspace{.1in}
Since we need the $L^2$ norm of $h$ on the left, we will choose $\alpha =\frac{2n-2}{n}$. Then

\vspace{.1in}

\noindent $(\displaystyle \int_{N_{\delta}}{h^2)}^{\frac{n-1}{n}}\leq c (\displaystyle \int_{N_{\delta}}|\nabla h^{\frac{2n-2}{n}}|+h^{\frac{2n-2}{n}}|{\mathcal H}_{\delta}|)\leq C \displaystyle \int_{N_{\delta}}|\nabla h||h|^{\frac{n-2}{n}}+c\displaystyle \int_{N_{\delta}}|h|^{\frac{2n-2}{n}}|{\mathcal H}_{\delta}|$

\vspace{.1in}

\hspace{.1 in} $\leq C (\displaystyle \int_{N_{\delta}}|\nabla h|^{p})^{1/p}(\displaystyle \int_{N_{\delta}}|h|^{\frac{nq-2q}{n}})^{1/q}+c(\displaystyle \int_{N_{\delta}}|h|^{\frac{2np'-2p'}{n}})^{1/p'}(\displaystyle \int_{N_{\delta}}|{\mathcal H}_{\delta}|^{q'})^{1/q'}$

\vspace{.1in}

Here $\frac {1}{p} +\frac {1}{q}= 1$ and $\frac {1}{p'}+\frac {1}{q'}=1$ and for our purposes we choose $p=\frac{2n}{n+2}$, $q=\frac{2n}{n-2}$, $p'=\frac{2n}{2n-2}$, and $q'=n$. Then the inequality becomes

\vspace{.1in}

$\leq C(\displaystyle \int_{N_{\delta}}|\nabla h|^{\frac{2n}{n+2}})^{\frac{n+2}{2n}}(\displaystyle \int_{N_{\delta}}|h^{2}|)^{\frac{n-2}{2n}}+c(\displaystyle \int_{N_{\delta}}|h^{2}|)^{\frac{n-1}{n}}(\displaystyle \int_{N_{\delta}}|{\mathcal H}_{\delta}|^{n})^{1/n})$

\vspace{.1in}

Since $||{\mathcal H}_{\delta}||_{L^n}$ is uniformly bounded, (here $n=3$) we can replace it with some small constant $C_1$. Then 

\vspace{.1in}

$(\displaystyle \int_{N_{\delta}}{|h|^2)}^{\frac{n-1}{n}}\leq  C(\displaystyle \int_{N_{\delta}}|\nabla h|^{\frac{2n}{n+2}})^{\frac{n+2}{2n}}(\displaystyle \int_{N_{\delta}}|h^{2}|)^{\frac{n-2}{2n}}+cC_1(\displaystyle \int_{N_{\delta}}|h^{2}|)^{\frac{n-1}{n}}$

\vspace{.1in}

\noindent which implies

\vspace{.1in}

$(1-cC_1)(\displaystyle \int_{N_{\delta}}{h^2)}^{\frac{n-1}{n}}\leq C(\displaystyle \int_{N_{\delta}}|\nabla h|^{\frac{2n}{n+2}})^{\frac{n+2}{2n}}(\displaystyle \int_{N_{\delta}}|h^{2}|)^{\frac{n-2}{2n}}$

\vspace{.1in}

$(\frac{1-cC_1}{C})(\displaystyle \int_{N_{\delta}}{|h^2|)}^{\frac{n-1}{n}-\frac{n-2}{2n} }\leq (\displaystyle \int_{N_{\delta}}|\nabla h|^{\frac{2n}{n+2}})^{\frac{n+2}{2n}}$

\vspace{.1in}

$(\frac{1-cC_1}{C})(\displaystyle \int_{N_{\delta}}{|h^2|)}^{1/2}\leq (\displaystyle \int_{N_{\delta}}|\nabla h|^{\frac{2n}{n+2}})^{\frac{n+2}{2n}}$

\vspace{.1in}

$(\frac{1-cC_1}{C})||h||_{L^2}\leq ||\nabla h||_{L^{(\frac{n+2}{2n})}}$

\vspace{.1in}

\noindent since $\frac{2n}{n+2}<2$, and $||\nabla h||_{L^{(\frac{n+2}{2n})}}\leq ||\nabla h||_{L^{2}} $ we finally get 

\vspace{.1in}

$(\frac{1-cC_1}{C})||h||_{L^2}\leq ||\nabla h||_{L^{2}}$ which proves the lemma.

\vspace{.1in}

Therefore the lemma 4.1 is proved.



{\lem : Let $f\in L^2_2(H_\delta)$ and $\psi \in L^2(H_\delta)$. Assume $\Delta_\delta f=\psi$ on $H_\delta$. Then we have 

\vspace{.1in}

\hspace{.2in} $||f||_{L^2_\delta}\leq C||\psi||_{L^2_\delta}$.}

\vspace{.1in}

{\bf Proof:} We multiply both sides of the equation $\Delta_\delta f=\psi$ by $f$ and integrate by parts, we get 




\vspace{.1in}

$\displaystyle\int_{B_\delta}{<d^*df,f>}$dvol$_\delta = \displaystyle\int_{B_\delta}{|df|^2}$dvol$_\delta = \displaystyle\int_{B_\delta}{<\psi ,f>}$dvol$_\delta$

\vspace{.1in}

$\leq (\displaystyle\int_{B_\delta}{|f|^2}$dvol$_\delta)^{1/2} (\displaystyle\int_{B_\delta}{|\psi |^2}$dvol$_\delta)^{1/2}$ 


\vspace{.1in}

Also by the eigenvalue estimate (Lemma 4.4) we get



\vspace{.1in}

$\displaystyle\int_{B_\delta}{|f|^2}$dvol$_\delta \leq C \displaystyle\int_{B_\delta}{|df|^2}$dvol$_\delta$

\vspace{.1in}

\noindent which implies $(\displaystyle\int_{B_\delta}{|f|^2}$dvol$_\delta)^{1/2}\leq C (\displaystyle\int_{B_\delta}{|\psi |^2}$dvol$_\delta)^{1/2}$ 


\vspace{.1in}

\noindent where $C$ is a constant that comes from the eigenvalue estimate which does not depend on $\delta$ and we get the lemma proved.

\vspace{.1in}

Next we will generalize the above estimate to general $p$.

{\lem : Let $f\in L^p_{2}(H_\delta)$ and $\psi \in L^p(H_\delta)$. Assume $\Delta_\delta f=\psi$ on $H_\delta$. Then we have 

\vspace{.1in}

\hspace{.2in} $||f||_{L^p_{\delta}}\leq C||\psi||_{L^p_{\delta}}$.}

\vspace{.1in}

{\bf Proof:} We multiply both sides of the equation $\Delta_\delta f=\psi$ by $f^{p-1}$ and integrate by parts, we get 

\vspace{.1in}

$\displaystyle\int_{B_\delta}{<d^*df,f^{p-1}>}dvol_\delta = \displaystyle\int_{B_\delta}{<df,d(f^{p-1})>}dvol_\delta = \displaystyle\int_{B_\delta}{<\psi ,f^{p-1}>}dvol_\delta$

\vspace{.1in}

$\leq (\displaystyle\int_{B_\delta}{|f|^p}dvol_\delta)^{\frac{p-1}{p}} (\displaystyle\int_{B_\delta}{|\psi |^q}dvol_\delta)^{\frac{1}{q}}$ 

\vspace{.1in}
\noindent where we have the identity $\frac{p-1}{p}+\frac{1}{q}=1$. We will also assume that $q\leq p$.


\vspace{.1in}

We also have 

\vspace{.1in}
$\displaystyle\int_{B_\delta}{<df,d(f^{p-1})>}$dvol$_\delta = \displaystyle\int_{B_\delta}{<df,(p-1)f^{p-2}df>}$dvol$_\delta$

\vspace{.1in}

\hspace{1.7in} $= \frac{4(p-1)}{p^2}\displaystyle\int_{B_\delta}{|d(f^{p/2})|^2}$dvol$_\delta$ 

\vspace{.1in}

Again by the eigenvalue estimate (Lemma 4.4) we get

\vspace{.1in}

$\frac{4(p-1)}{p^2}\displaystyle\int_{B_\delta}{|f^{p/2}|^2}$dvol$_\delta$$\leq C \frac{4(p-1)}{p^2}\displaystyle\int_{B_\delta}{|d(f^{p/2})|^2}$dvol$_\delta$  

\vspace{.1in}

\noindent which implies

\vspace{.1in}

$\displaystyle\int_{B_\delta}{|f^{p/2}|^2}$dvol$_\delta$$\leq C (\displaystyle\int_{B_\delta}{|f|^p}$dvol$_\delta)^{\frac{p-1}{p}} (\displaystyle\int_{B_\delta}{|\psi |^q}$dvol$_\delta)^{\frac{1}{q}}$ 

\vspace{.2in}

So we get $||f||_{L^p_\delta}\leq C||\psi||_{L^q_\delta}$. Since we chose $q\leq p$, this also implies 

\noindent $||f||_{L^p_\delta}\leq C||\psi||_{L^p_\delta}$ and again we get the lemma proved.

\vspace{.1in}

Next we will show that we have the basic elliptic estimates in the neck region. Note that the constant $C$ in these estimates depends on the geometry of the domain and as $\delta$ approaches to zero the neck region becomes almost singular. Therefore, we need to show that we can control that constant $C$ and obtain the elliptic estimates for the norms $||\cdot||_{L^2_{2,\delta}}$ and $||\cdot||_{L^2_{4,\delta}}$.

{\lem :  For $\Delta_{\delta}f=\psi$ there exists a constant $C$, independent of $\delta$ such that the following estimate holds:

\vspace{.1in} 

$||f||_{L^2_{2,\delta}}\leq C ||\psi||_{L^2_{\delta}}$  }

\vspace{.2in}

{\bf Proof:} First we write the Laplace operator in local coordinates, $\Delta_\delta f= \frac{1}{\sqrt { \tilde{g}}}\sum\limits_{j,k=1}^{n}\frac{\partial}{\partial x^j}(\sqrt { \tilde{g}}\tilde{g}^{jk}\frac{\partial}{\partial x^k})$. Here we will assume that $K$ is of dimension $n-2$. For simplicity we will also assume that $\mathcal A=1$ in the induced metric. Then we can write the metric in the neck region by 

\vspace{.1in}

$\tilde{g}_{\delta}|_{T_{\delta}\times K}=(1+\frac {\delta^4}{4(x^2+y^2)^2})(dxdx+dydy)+d{\kappa_1}d\kappa_1+...+d{\kappa_{n-2}}d\kappa_{n-2}$ 

\vspace{.1in}

\noindent and 

\vspace{.1in}

$\sqrt {det \tilde{g}_{\delta,ij}}=((1+\frac {\delta^4}{4(x^2+y^2)^2})(1+\frac {\delta^4}{4(x^2+y^2)^2}))^{1/2}=(1+\frac {\delta^4}{4(x^2+y^2)^2})$. 

\vspace{.1in}

For simplicity put $\mu=(1+\frac {\delta^4}{4(x^2+y^2)^2})$ . Using the metric we get $\Delta f=[ \mu^{-1}(\frac{\partial^2}{\partial x^2}+\frac{\partial^2}{\partial y^2}) + \frac{\partial^2}{\partial \kappa_1^2}+...+\frac{\partial^2}{\partial \kappa_{n-2}^2}]f = \psi$

\vspace{.1in}

Then we have
\vspace{.1in}

$(\frac{\partial^2 f}{\partial x^2}+\frac{\partial^2 f}{\partial y^2})+\mu(\frac{\partial^2f}{\partial \kappa_1^2}+...+\frac{\partial^2f}{\partial \kappa_{n-2}^2})=\overline\Delta f + (\frac {\delta^4}{4(x^2+y^2)^2})(\frac{\partial^2f}{\partial \kappa_1^2}+...+\frac{\partial^2f}{\partial \kappa_{n-2}^2})$

\vspace{.1in}

\hspace{2.27in} $=\mu \psi$

\vspace{.1in}
\noindent where $\overline\Delta= \frac{\partial^2}{\partial x^2} +\frac{\partial^2}{\partial y^2}+...+\frac{\partial^2}{\partial \kappa_{n-2}^2}$ is the standard n-dimensional Laplacian. 

\vspace{.1in}
Since $\overline\Delta$ has constant coefficients we have the following estimate for some $C$ independent of $\delta$ from Lemma 4.5 and basic elliptic estimates for flat Laplacian [2].

\vspace{.1in}

$||f||_{L_2^2}\leq C ||\mu \psi -  (\frac {\delta^4}{4(x^2+y^2)^2})(\frac{\partial^2f}{\partial \kappa_1^2}+...+\frac{\partial^2f}{\partial \kappa_{n-2}^2})||_{L^2}$

\vspace{.1in}

$\leq C||\mu \psi||_{L^2} + || (\frac {\delta^4}{4(x^2+y^2)^2})(\frac{\partial^2f}{\partial \kappa_1^2}+...+\frac{\partial^2f}{\partial \kappa_{n-2}^2})||_{L^2}$

\vspace{.1in}

$\leq C(($sup$|\mu|^2)^{1/2}\cdot ||\psi||_{L^2} + ($sup$|\frac {\delta^4}{4(x^2+y^2)^2}|^2)^{1/2}\cdot ||\frac{\partial^2f}{\partial \kappa_1^2}+...+\frac{\partial^2f}{\partial \kappa_{n-2}^2}||_{L^2})$

\vspace{.1in}

We have sup$|\mu|=1+\frac{\delta^4}{4\delta^4}=\frac{5}{4}$ and 
sup$|\frac {\delta^4}{4(x^2+y^2)^2}|=\frac{\delta^4}{4\delta^4}=\frac{1}{4}$ in the annulus defined by $\delta^2<(x^2+y^2)<\delta$.

\vspace{.1in}
Therefore

\vspace{.1in}
$||f||_{L_2^2}\leq C (C_1||\psi||_{L^2}+C_2 ||\frac{\partial^2f}{\partial \kappa_1^2}+...+\frac{\partial^2f}{\partial \kappa_{n-2}^2}||_{L^2})$

\vspace{.1in}
It remains to show that $||\frac{\partial^2f}{\partial \kappa_1^2}+...+\frac{\partial^2f}{\partial \kappa_{n-2}^2}||_{L^2}\leq C ||\psi||_{L_2}$. For a given test function $\tilde{v}\in C^1_0(U)$ on $U=U_\delta=\{\delta^2<x^2+y^2<\delta\}\times R^{n-2}\subset R^n$ we have 

\vspace{.1in}
${\displaystyle \int_U (\overline\Delta f)\tilde{v}+} \frac {\delta^4}{4(x^2+y^2)^2}(\frac{\partial^2f}{\partial \kappa_1^2}+...+\frac{\partial^2f}{\partial \kappa_{n-2}^2})\tilde{v} $dvol$_{\delta}= \displaystyle \int_U \mu \psi \tilde{v} $dvol$_{\delta}$

\vspace{.1in}
Since $\frac {\delta^4}{4(x^2+y^2)^2}$ is independent of $\kappa_1,...\kappa_{n-2}$ we can rewrite the above equation as follows:

\vspace{.1in}

${\displaystyle \int_U } f_x\tilde{v}_x+f_y\tilde{v}_y+\frac {\delta^4}{4(x^2+y^2)^2}(f_
{\kappa_1}\tilde{v}_
{\kappa_1}+...f_{\kappa_{n-2}}\tilde{v}_{\kappa_{n-2}})  $dvol$_\delta = -{\displaystyle \int_U }\mu \psi \tilde{v} $dvol$_{\delta}$

\vspace{.1in}
Fix $j\in \{1,..,n-2\}.$ Then

\vspace{.1in}

${\displaystyle \int_U}( f_{x\kappa_j}\tilde{v}_x+f_{y\kappa_j}\tilde{v}_y+(\frac {\delta^4}{4(x^2+y^2)^2})(f_
{\kappa_1 \kappa_j}\tilde{v}_
{\kappa_1}+...f_{\kappa_{n-2}\kappa_j}\tilde{v}_{\kappa_{n-2}}))$dvol$_\delta $

\vspace{.1in}

$=-{\displaystyle \int_U }(f_x\tilde{v}_{x\kappa_j}+f_y\tilde{v}_{y\kappa_j}+(\frac {\delta^4}{4(x^2+y^2)^2})(f_
{\kappa_1}\tilde{v}_
{\kappa_1\kappa_j}+...f_{\kappa_{n-2}}\tilde{v}_{\kappa_{n-2}\kappa_j}))  $dvol$_\delta $

\vspace{.1in}
$={\displaystyle \int_U}  \mu \psi \tilde{v}_{\kappa_j} $dvol$_{\delta}$

\vspace{.1in}

We will take $\tilde{v}=\nu ^2 (f_{\kappa_j})$ where $\nu$ is a cutoff function $\nu$ with uniformly bounded derivatives in $U$. By strict ellipticity and Cauchy Schwarz inequality we get

\vspace{.1in}

${\displaystyle \int_U} |f_{x\kappa_j}|^2+|f_{y\kappa_j}|^2+|f_{\kappa_1\kappa_j}|^2...+|f_{\kappa_{n-2}\kappa_j}|^2$dvol$_{\delta}$

\vspace{.1in}

$\leq C{\displaystyle \int_U}(f_{x\kappa_j}f_{x\kappa_j}\nu^2+ f_{y\kappa_j}f_{y\kappa_j}\nu^2+ \mu \nu^2 (f_{\kappa_1\kappa_j}f_{\kappa_1\kappa_j}+...+f_{\kappa_{n-2}\kappa_j}f_{\kappa_{n-2}\kappa_j})$dvol$_{\delta}$

\vspace{.1in}

\noindent where $C$ depends on the derivatives of $\nu$.

\vspace{.1in}

$\leq C{\displaystyle \int_U}  \mu \psi \nu^2 f_{\kappa_j\kappa_j}$dvol$_{\delta}\leq ||\mu \psi||_{L^2}\cdot ||f_{\kappa_j\kappa_j}||_{L^2}$ 

\vspace{.1in}

$\leq C||\psi||_{L^2}\cdot ||f_{\kappa_j\kappa_j}||_{L^2}$

\vspace{.1in}

Since $||f_{\kappa_j\kappa_j}||^2_{L^2}\leq {\displaystyle \int_U} (|f_{x\kappa_j}|^2+|f_{y\kappa_j}|^2+|f_{\kappa_1\kappa_j}|^2...+|f_{\kappa_{n-2}\kappa_j}|^2$dvol$_{\delta}$

\vspace{.1in}

\noindent we get $||f_{\kappa_j\kappa_j}||^2_{L^2}\leq C||\psi||_{L^2}\cdot ||f_{\kappa_j\kappa_j}||_{L^2}$ and hence $||f_{\kappa_j\kappa_j}||_{L^2}\leq C||\psi||_{L^2}.$

\vspace{.1in}

Finally $||\sum \limits_{j=1}^{k-2} f_{{k_j}{k_j}}||_{L^2}\leq \sum \limits_{j=1}^{k-2} ||f_{{k_j}{k_j}}||_{L^2}\leq C||\psi||_{L^2}$ and we get the result. Generalization this to $L^p$, for $p>2$ is straightforward.

\vspace{.1in}
Using a similar argument, one can improve the above inequality:

\vspace{.1in}

{\lem :  For $\Delta_{\delta}f=\psi$ there exists a constant $C$, independent of $\delta$ such that the following estimate holds:

\vspace{.1in} 

$||f||_{L^2_{4,\delta}}\leq C ||\psi||_{L^2_{2,\delta}}$  }

\vspace{.2in}

{\bf Proof:} As before we write the Laplace operator in local coordinates and we get
 \vspace{.1in}

$(\frac{\partial^2 f}{\partial x^2}+\frac{\partial^2 f}{\partial y^2})+\mu(\frac{\partial^2f}{\partial \kappa_1^2}+...+\frac{\partial^2f}{\partial \kappa_{n-2}^2})=\overline\Delta f + (\frac {\delta^4}{4(x^2+y^2)^2})(\frac{\partial^2f}{\partial \kappa_1^2}+...+\frac{\partial^2f}{\partial \kappa_{n-2}^2})$
\vspace{.1in}

\hspace{2.24in} $=\mu \psi$

\vspace{.1in}
\noindent where $\overline\Delta= \frac{\partial^2}{\partial x^2} +\frac{\partial^2}{\partial y^2}+...+\frac{\partial^2}{\partial \kappa_{n-2}^2}$ is the standard n-dimensional Laplacian with constant coefficients. 

\vspace{.1in}

Then we have the following estimates for flat Laplacian [2].

\vspace{.1in}

$||f||_{L_4^2}\leq C ||\mu \psi -  (\frac {\delta^4}{4(x^2+y^2)^2})(\frac{\partial^2f}{\partial \kappa_1^2}+...+\frac{\partial^2f}{\partial \kappa_{n-2}^2})||_{L_2^2}$

\vspace{.1in}

\hspace{.48in}$\leq C||\mu \psi||_{L_2^2} + || (\frac {\delta^4}{4(x^2+y^2)^2})(\frac{\partial^2f}{\partial \kappa_1^2}+...+\frac{\partial^2f}{\partial \kappa_{n-2}^2})||_{L_2^2}$

\vspace{.1in}


\vspace{.1in}
Since $||\nabla^2\mu||_{L^2}$, $||\nabla\mu||_{L^2}$ and $||\mu||_{L^2}$ are all uniformly bounded on $U=U_\delta=\{\delta^2<x^2+y^2<\delta\}\times R^{n-2}\subset R^n$ we get


\vspace{.1in}

$||f||_{L_2^4}\leq C (C_1||\psi||_{L_2^2}+C_2 ||\frac{\partial^2f}{\partial \kappa_1^2}+...+\frac{\partial^2f}{\partial \kappa_{n-2}^2}||_{L_2^2})$

\vspace{.1in}
As before it remains to show that $||\frac{\partial^2f}{\partial \kappa_1^2}+...+\frac{\partial^2f}{\partial \kappa_{n-2}^2}||_{L_2^2}\leq C ||\psi||_{L_2^2}$. By Lemma 4.7 we have $||\frac{\partial^2f}{\partial \kappa_1^2}+...+\frac{\partial^2f}{\partial \kappa_{n-2}^2}||_{L^2}\leq C ||\psi||_{L^2}$. So we only need to show that $||\nabla (\frac{\partial^2f}{\partial \kappa_1^2}+...+\frac{\partial^2f}{\partial \kappa_{n-2}^2})||_{L^2}\leq C ||\nabla \psi||_{L^2}$ and similarly $||\nabla^2 (\frac{\partial^2f}{\partial \kappa_1^2}+...+\frac{\partial^2f}{\partial \kappa_{n-2}^2})||_{L^2}\leq C ||\nabla^2 \psi||_{L^2}$. 

\vspace{.1in}
In order to prove the first inequality we will follow the same arguments as before. The only change will be to replace the test function $\tilde{v}$ with $\nu^2 ({\nabla}^* \nabla f_{\kappa_j})$  where $\nu$ is a cutoff function $\nu$ with uniformly bounded derivatives in $U$.














\vspace{.1in}

${\displaystyle \int_U} |\nabla f_{x\kappa_j}|^2+|\nabla f_{y\kappa_j}|^2+|\nabla f_{\kappa_1\kappa_j}|^2...+|\nabla f_{\kappa_{n-2}\kappa_j}|^2$dvol$_{\delta}$

\vspace{.1in}

$\leq C{\displaystyle \int_U}(\nabla f_{x\kappa_j}\nabla f_{x\kappa_j}\nu^2+ \nabla f_{y\kappa_j}\nabla f_{y\kappa_j}\nu^2$

\hspace{1in} $+\mu \nu^2 (\nabla f_{\kappa_1\kappa_j}\nabla f_{\kappa_1\kappa_j}+...+\nabla f_{\kappa_{n-2}\kappa_j}\nabla f_{\kappa_{n-2}\kappa_j})$dvol$_{\delta}$

\vspace{.1in}


$\leq C{\displaystyle \int_U}  \mu \psi \nu^2 {\nabla}^* \nabla f_{\kappa_j\kappa_j}$dvol$_{\delta}\leq ||\mu \psi||_{L^2}\cdot ||\nabla f_{\kappa_j\kappa_j}||_{L^2}$ 

\vspace{.1in}

$\leq C||\nabla \psi||_{L^2}\cdot ||\nabla f_{\kappa_j\kappa_j}||_{L^2}$

\vspace{.1in}

Since $||\nabla f_{\kappa_j\kappa_j}||^2_{L^2}\leq {\displaystyle \int_U} (|\nabla f_{x\kappa_j}|^2+|\nabla f_{y\kappa_j}|^2+|\nabla f_{\kappa_1\kappa_j}|^2...+|\nabla f_{\kappa_{n-2}\kappa_j}|^2$dvol$_{\delta}$

\vspace{.1in}

\noindent we get 

\vspace{.1in}

\noindent $||\nabla f_{\kappa_j\kappa_j}||^2_{L^2}\leq C||\nabla \psi||_{L^2}\cdot ||\nabla f_{\kappa_j\kappa_j}||_{L^2}$ and so $||\nabla f_{\kappa_j\kappa_j}||_{L^2}\leq C||\nabla \psi||_{L^2}.$

\vspace{.1in}

Since $||\sum \limits_{j=1}^{k-2} \nabla f_{{k_j}{k_j}}||_{L_2}\leq \sum \limits_{j=1}^{k-2} ||\nabla f_{{k_j}{k_j}}||_{L_2}\leq C||\nabla \psi||_{L_2}$ we get the inequality $||\nabla (\frac{\partial^2f}{\partial \kappa_1^2}+...+\frac{\partial^2f}{\partial \kappa_{n-2}^2})||_{L^2}\leq C ||\nabla \psi||_{L^2}$.

\vspace{.1in}

The proof of the inequality $||\nabla^2 (\frac{\partial^2f}{\partial \kappa_1^2}+...+\frac{\partial^2f}{\partial \kappa_{n-2}^2})||_{L^2}\leq C ||\nabla^2 \psi||_{L^2}$ follows exactly the same way. The combination of these inequalities will prove lemma 4.8.

\vspace{.1in}

Using Lemmas 4.7 and 4.8 we can finally state the following theorem about the uniform invertibility of the Laplacian.

\vspace{.1in}

{\thm : The linear operator $d^*_{\delta}d=\Delta_{\delta} :C^\infty(H_{\delta})\rightarrow C^\infty(H_{\delta})$ has a right inverse $Q_{\delta}$ under the norm $||.||_{\delta}$ in the sense that there exists a constant $C$, independent of $\delta$, such that for sufficiently small $\delta$,

\vspace{.1in}

$||Q_{\delta}g||_{\delta,k+2,p}<C||g||_{\delta,k,p}$ for any $g\in L^p(H_{\delta})$ and $k\leq 2$.}

\vspace{.1in}

\vspace{.1in}

Note that we assume in the beginning that the almost special Lagrangian is deformed under Hamiltonian deformations, in other words for any given $\eta$ there exists a function $g$ satisfying $dg=\eta$. Therefore in this case the invertibility of Laplacian $\triangle_{\delta}$ on smooth functions will be equivalent to the invertibility of $d^*_{\delta}$ on (exact) one forms. So we can restate the above theorem as follows:
  
\vspace{.1in}

{\thm : The linear operator $*d^*_{\delta} :\Gamma(N(H_{\delta}))\rightarrow \Omega^3(H_{\delta})$ has a right inverse $P_{\delta}$ under the norm $||.||_{\delta}$ in the sense that there exists a constant $C$, independent of $\delta$, such that for sufficiently small $\delta$,

\vspace{.1in}

$||P_{\delta}\eta||_{\delta,k+1,p}<C||\eta||_{\delta,k,p}$ for any $\eta\in L^p(\Omega^3(H_{\delta}))$ and $k\leq 2$.}

\vspace{.1in}

\section{Implicit Function Theorem}

\vspace{.1in}

In this section we will write the Taylor expansion for the special Lagrangian equation:
$\mathcal F(\exp_{H_\delta}w)=\mathcal F(H_\delta )+D_{\delta}\mathcal F(w)+\mathcal O(w)$ and show that the non-linear term $\mathcal O(w)$ satisfies the contraction mapping principle. Then by the Implicit Function Theorem, we show that there exists a unique $w$ that satisfies the equation $\mathcal F(H_\delta )+D_{\delta}\mathcal F(w)+\mathcal O(w)=0$ and therefore there exists a real smooth special Lagrangian submanifold near $L$.

\vspace{.1in}

The map $\mathcal F$ has the following Taylor expansion:

\vspace{.1in}

\hspace{.1in} $\mathcal F(\exp_{H_\delta} w)=\mathcal F(H_\delta )+D_{H_\delta}\mathcal  F(w)+\mathcal O(w)$

\vspace{.1in}

\noindent where $D_{H_\delta} \mathcal F(w)$ is the linearization of the deformation map $\mathcal F$ and $\mathcal O(w)$ is the non-linear term.

\vspace{.1in}

It is obvious that for given $\delta$ the solution set of the following equation is the set of vector fields that corresponds to the real special Lagrangian submanifold near $H_\delta$.

\vspace{.1in}

\hspace{.1in} $\mathcal F(\exp_{H_\delta} w)=\mathcal F(H_\delta )+D_{H_\delta} \mathcal F(w)+\mathcal O(w)=0$

\vspace{.1in}

As before we can replace the one form $w=\nabla h$ with a smooth function and rewrite the equation in terms of $h$ as follows: 

\vspace{.1in}

\hspace{.1in} $\mathcal F(\exp_{H_\delta} h)=\mathcal F(H_\delta )+D_{H_\delta}\mathcal  F(h)+\mathcal O(h)=0$

\vspace{.1in}

By Theorem 4.9, the operator $D_{H_\delta}\mathcal F$ has a well defined inverse $G_\delta$. If we apply $G_\delta$ to the equation above we get 

\vspace{.1in}

$G_\delta \mathcal F(H_\delta )+h+G_\delta \mathcal O(h)=0$ which implies $-G_\delta \mathcal F(H_\delta )-G_\delta \mathcal O(h)=h$.

\vspace{.1in}

For simplicity, we will call the operator $-G_\delta \mathcal F(H_\delta )-G_\delta \mathcal O(h)=\mathcal W_\delta  (h)$. Then if $h$ is a solution to $\mathcal W_\delta h=h$, it is a solution to the equation $\mathcal F(h)=0$.
It is a basic fact that an operator $\mathcal W_\delta $ has a fixed point if the contraction mapping principle holds for $\mathcal W_\delta $.

\vspace{.1in}

First we will recall a basic fact in Banach space calculus.

{\lem : Let $\mathcal F:X'\rightarrow Y'$ be a smooth map between Banach spaces with the Taylor expansion

\vspace{.1in}

\noindent $\mathcal F(\exp_H x)=\mathcal F(H)+D_H\mathcal F(x)+\mathcal O(x)=0$  with dim(ker $D\mathcal F(0))< \infty$.

\vspace{.1in}

Also assume that $G:Y'\rightarrow X'$ to be a right inverse  $D\mathcal F(0)\circ G=id_B$. And let $C>0$ be a constant such that the nonlinear terms $\mathcal O(x),\mathcal O(y)$ satisfy the contraction property

\vspace{.1in} 

\noindent $||G\mathcal O(x)-G\mathcal O(y)||\leq C(||x||+||y||)||x-y||$ for all $x,y \in B_\epsilon(0)$ with $\epsilon=\frac{1}{8C}$.

\vspace{.1in}

Then $||G\mathcal F(0)|| \leq \epsilon /2$ will imply the unique existence of a zero, $x_0$ of $\mathcal F$ which satisfies $||x_0||\leq 2||G\mathcal F(0)||$. }

\vspace{.1in}

{\bf Proof:} The proof is an application of the contraction mapping principle. We will consider the map

\vspace{.1in}

$\mathcal W:A\rightarrow A$, $\mathcal W(x)=-G(\mathcal F(0)+\mathcal O(x)),$ as $x\in B_\epsilon(0)$ we get the following estimate

\vspace{.1in}

$||\mathcal W(x)||\leq \epsilon /2 +C{\epsilon}^2=\epsilon(1/2+1/8)<\epsilon$

\vspace{.1in}

\noindent and since 

\vspace{.1in}

$||\mathcal W(x)-\mathcal W(y)||=||G\mathcal O(x)-G\mathcal O(y)||\leq C2\epsilon||x-y||$

\vspace{.1in}

\hspace{1.04in} $=\frac{2}{8}||x-y||\leq \frac{1}{2}||x-y||$

\vspace{.1in}

\noindent we can conclude that $\mathcal W(x)$ is a contraction on $B_\epsilon(0)$. Hence we find a unique $x_0\in B_\epsilon(0)$ satisfying $\mathcal W(x_0)=x_0$.

\vspace{.1in}

Moreover we obtain the following estimate

\vspace{.1in}

$||x_0||-||\mathcal W(0)||\leq ||\mathcal W(x_0)-\mathcal W(0)|| \leq 1/2||x_0||$ and deduce the inequality $||x_0||\leq 2||G\mathcal F(0)||$. Hence we proved Lemma 5.1.


\vspace{.2in}

The next step will be to show that the operator $\mathcal W_\delta (h)=-G_\delta \mathcal F(H)- G_\delta \mathcal O(h)$ satisfies the following inequality:

\vspace{.1in}

$||\mathcal W_\delta f-\mathcal W_\delta g||_{L^2_{2,\delta}}\leq k.||f-g||_{L^2_{2,\delta}}$ for some $k<1$.

\vspace{.1in}

\noindent where $||\mathcal W_\delta f-\mathcal W_\delta g||_{L^2_{2,\delta}}=||-G_\delta \mathcal F(H)-G_\delta \mathcal O(f)+G_\delta \mathcal F(H)-G_\delta \mathcal O(g)||_{L^2_{2,\delta}}$

\vspace{.1in}

\hspace{1.4in}$=||-G_\delta \mathcal O(f)+G_\delta \mathcal O(g)||_{L^2_{2,\delta}}$

\vspace{.1in}

We will do this by first showing that
 
\vspace{.1in}

$||\mathcal O(f)-\mathcal O(g)||_{L^2_{2,\delta}}=
||\frac{{\partial}^2f}{\partial x \partial x}\frac{{\partial}^2f}{\partial y \partial y}-\frac{{\partial}^2f}{\partial x \partial y}\frac{{\partial}^2f}{\partial y \partial x}-\frac{{\partial}^2g}{\partial x \partial x}\frac{{\partial}^2g}{\partial y \partial y}+\frac{{\partial}^2g}{\partial x \partial y}\frac{{\partial}^2g}{\partial y \partial x}||_{L^2_{2,\delta}}$

\vspace{.1in}

\hspace{1.3in} $\leq C||f-g||_{L^2_{4,\delta}}$

\vspace{.1in}

We will prove this inequality term by term. For simplicity we will write $\frac{{\partial}^2f}{\partial x \partial x}$ as $f_{ xx}$.

\vspace{.1in}

$||f_{ xx}f_{ yy}-g_{xx}g_{ yy}||_{L^2_{2,\delta}}=||f_{ xx}f_{ yy}-f_{ xx}g_{ yy}+f_{ xx}g_{ yy}-g_{ xx}g_{ yy}||_{L^2_{2,\delta}}$

\vspace{.1in}

\noindent $\leq C(||f_{xx}f_{ yy}-f_{ xx}g_{ yy}||_{L^2_{2,\delta}}
+||f_{xx}g_{ yy}-g_{ xx}g_{ yy}||_{L^2_{2,\delta}})$  (triangular inequality)

\vspace{.1in}

\noindent $\leq C(||f_{ xx}||_{L^2_{2,\delta}}||f_{ yy}-g_{yy}||_{L^2_{2,\delta}}
+ ||g_{ yy}||_{L^2_{2,\delta}}||f_{ xx}-g_{ xx}||_{L^2_{2,\delta}})$

\vspace{.1in}

\noindent $=C(||f_{ xx}||_{L^2_{2,\delta}}||(f-g)_{yy}||_{L^2_{2,\delta}}+||g_{ yy}||_{L^2_{2,\delta}}||(f-g)_{ xx}||_{L^2_{2,\delta}})$

\vspace{.1in}

\noindent $\leq C(||f||_{L^2_{4,\delta}}.||f-g||_{L^2_{4,\delta}}+||g||_{L^2_{4,\delta}}||f-g||_{L^2_{4,\delta}})$ (Elliptic estimates)

\vspace{.1in}

\noindent $=C(||f||_{L^2_{4,\delta}}+||g||_{L^2_{4,\delta}})||f-g||_{L^2_{4,\delta}}$

\vspace{.1in}

\noindent and similarly we get

\vspace{.1in}

\noindent $||f_{xy}f_{ yx}-g_{ xy}g_{yx}||_{L^2_{2,\delta}}\leq C(||f||_{L^2_{4,\delta}}+||g||_{L^2_{4,\delta}})||f-g||_{L^2_{4,\delta}}$

\vspace{.1in}

\noindent and we combine all these terms to get 

\vspace{.1in}

\noindent $||-\mathcal O(g)+ \mathcal O(f)||_{L^2_{2,\delta}}\leq C(||f||_{L^2_{4,\delta}}+||g||_{L^2_{4,\delta}})||f-g||_{L^2_{4,\delta}}$

{\lem : There exists a constant $\epsilon>0$ independent of the gluing parameter $\delta$ such that for sufficiently small $\delta$ and for any $\gamma=-G_\delta \mathcal F(H_\delta) \in  L^2_4(H_\delta)$ the equation:

\vspace{.1in}

\hspace{.1in} $h+\mathcal O_\delta(h)=\gamma $

\vspace{.1in}
\noindent has a unique small solution $h\in L^2_4(H_\delta)$ in the ball $||h||_{L^2_{2,\delta}}\leq \epsilon/2$. In other words $\mathcal F(\exp_{H_\delta} h)$ is a special Lagrangian submanifold. Moreover $h$ satisfies $||h||_{L^2_{2,\delta}}\leq 2||\gamma ||_{L^2_{2,\delta}}$.}

\vspace{.1in}

Here $\gamma $ represents the error term.

\vspace{.1in}

{\bf Proof:} We will apply the contraction principle to the operator $\mathcal W_\delta (h)=-G_\delta \mathcal F(H)- G_\delta \mathcal O(h)$ defined in a small ball around 0 in $L^2_4(H_\delta)$.

\vspace{.1in}

Our estimate above and Lemma 4.8 implies that 

\vspace{.1in}
\noindent $||\mathcal W(h_1)-\mathcal W(h_2)||_{L^2_{2,\delta}}\leq C(||G(h_1)||_{L^2_{4,\delta}}+||G(h_2)||_{L^2_{4,\delta}})||G(h_1-h_2)||_{L^2_{4,\delta}}$

\vspace{.1in}

\hspace{1.28in} $\leq C(||h_1||+||h_2||)_{L^2_{2,\delta}}||(h_1-h_2)||_{L^2_{2,\delta}}$.

\vspace{.1in}

If we choose $\epsilon<C/8$ for $||\mathcal W_\delta(0)||_{L^2_{2,\delta}}\leq \epsilon/2$ then $\mathcal W_\delta$ will be a contraction on $B_\epsilon (0)$. Therefore there exists a unique fixed point $h$ of $\mathcal W_\delta$ satisfying $||h||_{L^2_{2,\delta}}\leq 2||W_\delta(0)||_{L^2_{2,\delta}}$.



{\rem :} Note that when we constructed the almost special Lagrangian we showed that it could be represented as a graph of some smooth function. 
Then among these graphs, the solution set of an elliptic equation gives the real special Lagrangian submanifolds. This implies that the special Lagrangian submanifold that we obtain as the fixed point of the map $\mathcal W$ is in fact smooth and it can not be the same submanifold as the original singular submanifold $L$. 

\vspace{.1in}

Finally we can conclude the following theorem for the singularity of the form $z_1.\overline {z_2}=0$.

\vspace{.1in}

{\thm: Given a connected immersed special Lagrangian submanifold $L^3$ of a Calabi-Yau manifold $X^6$ with a particular irreducible self intersection $K$ of codimension-two (singularity of type $z_1.\overline {z_2}=0$) it can be approximated by a sequence of smooth special Lagrangian submanifolds and therefore $L$ is a limit point in the moduli space. }

\vspace{.1in}

{\rem :} We hope to extend this result to the cusp-type singularities in a Calabi-Yau manifold. The details of this construction for both the Calabi-Yau and the symplectic case will appear somewhere else.

\vspace{.1in}

\small
{\em Acknowledgements.} 


Special thanks are to Gang Tian for all his guidance and encouragement during the course of this work. The author would also like to thank Allen Back, Yuri Berest, Laurent Saloff-Coste, Jos{\'e} Escobar and Peng Lu for many useful conversations. Also thanks to Dominic Joyce, Tom Parker and Zhengfang Zhou  for their useful comments in the earlier versions of this paper.

\vspace{.4in}


\begin{thebibliography}{[FP]}

\bibitem[1]{bf} Bryant, R.L. {\em Some examples of special Lagrangian Tori}, math.DG/9902076

\bibitem[2]{bf} Gilbarg, D. and Trudinger S.N.{\em Elliptic Partial Differential Equations of Second Order}, Springer-Verlag (New York 1977)


\bibitem[3]{bf} Gross, M. {\em Special Lagrangian Fibrations I: Topology}, Integrable Systems and Algebraic Geometry (Kobe/Kyoto 1997), 156-193 World Scientific.



\bibitem[4]{bf} Gross, M. {\em Special Lagrangian Fibrations II: Geometry}, alg-geom/9809072


\bibitem[5]{ab} Hitchin, N. {\em The moduli space of special Lagrangian submanifolds}, dg-ga/9711002



\bibitem[6]{b} Harvey, F.R. and Lawson, H.B. {\em Calibrated Geometries}, Acta. Math. {\bf 148} (1982), 47-157

\bibitem[7]{ab} Ionel, E.N. and Parker, T.H. {\em The Symplectic Sum Formula for Gromov-Witten Invariants}, dg-ga/0010217


\bibitem[8]{b} Joyce, D. {\em On counting special Lagrangian homology 3-spheres}, hep-th/9907013


\bibitem[9]{b} Joyce, D. private communication.


\bibitem[10]{b} Liu, G. {\em Associativity of Quantum Multiplication}, Commun.Math.Phys. {\bf 191} (1998), 265-282



\bibitem[11]{ab} Lu, P. {\em Special Lagrangian Tori on a Borcea-Voisin Threefold}, dg-ga/9902063



\bibitem[12]{ab} McDuff, D. and Salamon, D.A. {\em J-holomorphic Curves and Quantum Cohomology}, University Lecture Series {\bf 6}, American Mathematical Society, Providence, RI.



\bibitem[13]{bf} McLean, R.C. {\em Deformations of calibrated submanifolds}, Comm. Anal. Geom. {\bf 6} (1998), 705-747 



\bibitem[14]{bf} Mazzeo, R., Pollack, D., Uhlenbeck, K., {\em Connected Sum Constructions For Constant Scalar Curvature Metrics}, Topological Methods in Nonlinear Analysis, Journal of the Juliusz Schauder Center, {\bf Vol.6} (1995), 207-233



\bibitem[15]{bf} Ruan, Y. and Tian, G. {\em A Mathematical Theory of Quantum Cohomology}, J.Differential Geometry {\bf 42} (1995), 259-367



\bibitem[16]{bf} Salur, S. {\em Deformations of Special Lagrangian Submanifolds}, Communications in Contemporary Mathematics {\bf Vol.2, No.3} (2000), 365-372 

\bibitem[17]{bf} Simon, L. {\em Lectures on Geometric Measure Theory}, Proceedings of The Center for Mathematical Analysis Australian National University {\bf Vol.3} (1983)

\bibitem[18]{bf} Strominger, A., Yau, S.T. and Zaslow, E., {\em Mirror Symmetry is T-Duality}, Nucl. Phys. {\bf B479} (1996), 243-259



\end{thebibliography}
\end{document}